\documentclass[12pt]{article}

\usepackage{amsmath}

\usepackage{amssymb}

\usepackage{exscale}

\usepackage{amsthm}

\usepackage{epsfig}

\usepackage{verbatim}

\theoremstyle{plain}

\newcommand{\R}{{\mathbb{R}}}

\begin{document}

\title{Theory of zones on Zeeman manifolds\\ A new
approach to the infinities of QED}

\date{}

\author{Zolt\'an Imre Szab\'o \thanks{Lehman College of CUNY, Bronx, 
NY, 10468, and R\'enyi Institute, Budapest, Hungary.}}

\maketitle

\begin{abstract} The Zeeman-Hamilton operator of free 
charged particles are identified with the Laplacians of certain 
Riemannian manifolds, called
Zeeman manifolds. The quantum Hilbert space, $\mathcal{H}$, decomposes
into subspaces (Zeeman zones) which are invariant under the action
both of the Zeeman operator and the natural Heisenberg group 
representation. Thus a well defined particle 
theory and zonal geometry can be developed on each zone separately. 
The most surprising result
is that quantities those divergent
on the global setting appear to be finite on the
zonal setting. Even the zonal Feynman integral is well defined.

This zonal interpretation of particles has fundamental 
effect both on the physical and mathematical view of these objects.
The points are non-existing on a zone, for instance. 
One should introduce the concept of zonal point-spread, 
defined by certain wave functions. 
As a result, the zonal particles are not point- but 
point-spread-objects. The theory 
developed below includes explicit computation of objects
such as the waves defining the point-spreads, the zonal Wiener-Kac 
and Dirac-Feynman 
flows (which then define the corresponding measures on the
path-spaces), and the corresponding zonal Feynman-Kac formulas. 
It will be
tested against several well known effect such as the Aharanov-Bohm 
effect and Lamb shift. There is also
explained why these extended charged particles do not blow up.
 
\end{abstract}

\maketitle

\section{Introduction}

The confusing infinities (divergent integrals), stubbornly 
present in calculations since the very beginning of quantum
electrodinamics, are controlled by {\it renormalization} today. 
This pertubative tool provides
the desired finite quantities by differences of
infinities. 
The legacy of difficulties came from concepts such as {\it point mass}
and {\it point charge} of classical electron
theory, which provided the first warning that a point 
electron will have infinite electromagnetic
self-mass; the mass $e^2/6\pi ac^2$ for a surface distribution
of charge with radius $a$ blows up for $a\to 0$. In quantum
field theory the Hamiltonian of the field is proportional to this 
electromagnetic self-mass. This is why this infinity launched one 
of the deepest crisises in the history of physics.

Mathematically speaking, the infinities mostly due to the infinite
trace of kernels such as the Wiener-Kac kernel $e^{-tH}$
(providing the fundamental solution of the heat equation),
or, Dirac-Feynman kernel $e^{-tH\bold i}$ (providing the fundamental
solution of the Schr\"odinger equation). As a result, 
they assign infinite measures to 
physical objects such as
self-mass, self-charge, e.t.c.. 
Also the Feynman measure, which is analogous 
to the well defined Wiener-Kac measure on the path-spaces, 
requires renormalization.

This paper gives a new non-perturbative approach to this problem.
In the first step the quantum Hilbert space $\mathcal H$ 
(on which the quantum
Hamilton operator $H$ is acting) is decomposed into the direct
sum of $H$-invariant subspaces called
Zeeman zones. Then, all the actions such as the heat- or Feynman-flow 
are considered on these invariant subspaces separately. 
It turns out that both the Wiener-Kac and Dirac-Feynman kernels are 
of the trace class on each zone, furthermore, both define the 
corresponding zonal measures on the path-spaces rigorously.

Three types of Hamilton operators are considered in this paper.
(1) The classical
Zeeman operator $H_Z$ corresponding to free charged particles with
exterior (orbiting) spin. This spin is exhibited 
by the angular momentum operator 
associated with the magnetic dipole moment. (2) The Pauli operator, 
which provides also inner spin to the orbiting spin. (3) The 
relativistic Dirac operator.
These Hamiltonians are the most important ones in quantum physics. 
They were primarily introduced for explaining the Zeeman effect. 

There are several new features also in the manner by which these 
operators are introduced in this paper. The original historic model 
for establishing the 
Zeeman operator was a charged particle orbiting in a constant
magnetic field which is perpendicular to the plane where the
particle is orbiting. The Hamilton function of this 
system can be found by the Maxwell equations. 
In the hypothesis of particles having orbiting spin one assumes that
the Hamiltonian of charged particles are exactly of this form. This
hypothesis, so to speak, adds a little magnet to the particle.
Assumption $E=0$, imposed for the constant electromagnetic field, fixes
an inertial system in which the orbiting particle is considered. 
In fact, in all other systems the field appears with a non-vanishing 
electric field $E$. By this observation, all considerations developed
here become relativistic. Even the Pauli operator, operating on
spinors having two components, will be identified with the Dirac 
operator acting on four component spinors. This is why they are
called Pauli-Dirac ($\mathcal{PD}$) operators. 

An other new interpretation is that the Zeeman Hamiltonian 
operator, corresponding to the above Hamiltonian function, is 
identified with the Laplacian of a non-compact Riemannian manifold. 
In case of a single particle orbiting in the $(x,y)$-plane this manifold
is nothing but the Heisenberg group, parametrized by $(x,y,t)$, whose
center, parametrized by $t$, is periodic established by 
factorization with a lattice. The Riemannian metric, $g$, is the 
natural left-invariant metric on this manifold. Despite the
positive definite metric, the considerations are still relativistic, 
due to the above observation. 

One can generalize this idea
to higher dimensional Heisenberg groups. This model corresponds to more 
particles orbiting in the same inertial system determined by the same
constant magnetic field. Much more interesting is the generalization
to two-step nilpotent Lie groups whose center is factorized by a
co-compact lattice. This model corresponds to more particles being in
different constant magnetic fields, thus each of them determines
an individual inertial system with individual time. 
This mathematical model matches
Dirac's multi-time theory. All these Riemannian manifolds are called
Zeeman manifolds.

It is well known that quantum physics is lacking of a
comprehensive mathematical model
which is able to explain the big diversity of phenomenons experienced
in quantum theory. (The Schr\"odinger equation itself is not enough
for explaining many of these experiences.) The
above Riemannian manifolds partially furnish this missing piece, 
serving as a
relativistic space-time background for the theory presented here.
It should be pointed out that the fixed inertial system appearing
in the model resolves the concern of Pauli, brought up against Dirac's
relativistic electron theory. Indeed,
beside the hole theory, Pauli strongly criticized also 
the probabilistic argument developed on the relativistic 
space-time.
According to Pauli, probabilistic amplitudes do not have any meaning
on the Minkowski space-time, since they must be defined on the space. 
Note that the above depicted space-time concept really offers a 
resolution for this concern. In fact, 
the probabilistic theory will be developed on the space determined
by the fixed inertial system. The constant magnetic field, 
so to speak, serves as a
bridge for descending from the Minkowski 4-space to the 3-space, needed
for the probabilistic interpretations.

It should be pointed out that the
mathematical model depicted so far does not go beyond the 
Schr\"odinger equation,
therefore, it is not enough for explaining the effects which are not
in the scope of the classical theory.
The lacking piece what should be furnished yet 
is the spectral Zeeman zone decomposition, $\mathcal H=
\oplus_{a=0}^\infty \mathcal H^{(a)}$, of the quantum Hilbert space. 
On the complex $z=(x,y)$-plane the total quantum Hilbert space, 
$\mathcal H$, is the $L^2$-Hilbert space
of complex valued functions defined on $\mathbb C$, where the
inner product is defined by a Gaussian density 
$e^{-\lambda z \overline z}$. This Hilbert space is spanned by the
polynomials written in terms of $z$ and $\overline z$. 

The zones are
described in many different levels in this paper. They can be defined
by Gram-Schmidt orthogonalization such that the starting zone, 
$\mathcal H^{(0)}$ 
is the holomorphic zone spanned by the holomorphic
polynomials. Then 
$\mathcal H^{(1)}$
is defined by applying the Gram-Schmidt orthogonalization to the 
function space 
$G^{(1)}$ spanned by functions of the form 
$\overline zh(z)$, where $h(z)$ is an arbitrary holomorphic polynomial.
The procedure is continued by successive orthogonalization of 
spaces
$G^{(a)}$ spanned by functions of the form 
$\overline z^{a}h(z)$.
 
There is shown in the paper that the zones on the complex plane
are the irreducible invariant subspaces of the natural Heisenberg
group representation, thus each zone is suitable for establishing its
own physics. In higher dimensional cases, describing multiple of 
particles, the irreducibility fails, however, in this case a zone
further decomposes into subspaces by which exclusion principles and 
particles such as Bosonic and Fermionic ones can be defined.

Yet an other description explores the fact that the zones are invariant
under the action of the Hamilton operator and, by using explicit 
spectrum computations, the zones are introduced by the spectrum of 
the magnetic dipole moment operator. This technique actually uses
polarization. A zone corresponds to a magnetic state of 
the particles.

All the important spectral theoretical objects on Zeeman manifolds are
explicitly established. They include the spectrum,
the zonal projection operators, the zonal Wiener-Kac and
Dirac-Feynman kernels, and the zonal partition functions. 
Also the zonal Feynman-Kac type
formulas, for both the Wiener-Kac and
Feynman measure, are established in a novel form. One of the 
spectaculars of this theory is that all the quantities which appear
as infinities regarding the total space $\mathcal H$ become finite ones
on the zonal setting. Particularly, also the zonal Feynman measures are
well defined which can be explicitly computed.

The consideration of particles as zonal objects has a deep effect
also on the mathematical approach to the geometry formed within a zone.
The most important impact is that there do not exist points in the
original sense on a zone. By using the total Hilbert space $\mathcal H$
the points, $X$, can be identified with the Dirac delta distribution
$\delta_X(Y)=\sum \phi_i(X)\overline{\phi}(Y)$, where $\{\phi_i\}$ is
an orthonormal basis on $\mathcal H$. Regarding a zone, 
$\mathcal H^{(a)}$, this distribution can be defined by those functions,
$\phi_i^{(a)}$, which are in the zone. In other words, on a zone the 
points can be introduced by the kernel function belonging to the 
projection regarding the zone.
On the holomorphic zone, for instance, this kernel is the well known
Bergman kernel $\delta_x^{(0)}(y)=qe^{x\overline y-(1/2)(x^2+y^2)}$. It
is fascinating to see that one gets the kernel of 
projection onto an other zone, $\mathcal H^{(a)}$, just 
by multiplying the Bergman kernel by the Laguerre polynomial 
$L_a(|x-y|^2)$. Thus the zonal points are point-spreads which are
explicitly described by the above wave functions.

In physics the idea of matter-wave was 
developed by de Broglie, whose theory got its final form in the
Schr\"odinger equation. However, the Schr\"odinger theory still
needs the point and other relating concepts such as point-mass
and point-charge (see Weisskopf's explanation in the next section),
which then lead to the confusing infinities in quantum field theory.
Mathematics did not follow de Broglie's idea and it
still enforced the point concept,
from which the de Broglie theory departed.
The contradiction between the two approach
had been lifted by the duality principle, asserting that matter
manifests itself sometime as wave while in other cases as point 
particle.

The zonal theory, called {\it point spread geometry}, establishes
de Broglie's idea also on mathematical level. The point spreads on a 
zone are considered to be the most compressed mathematical objects which
still spread out over the whole space. They are the most compressed
states possible for a zonal particle.
These physical objects can be also in other states described by the
zonal wave functions. The question remains if a charge spread is stable
if it spreads as a point spread. It is pointed out that these charged 
point-spreads are not stable (they blow up), however, the zonal
Dirac-Feynman flow corresponding to the zonal Schr\"odinger equation
moves these spreads to stable, so called, solid charge spreads. The
duality principle is overtaken also to the zonal theory. The zonal 
objects manifests themselves sometimes by the zonal wave functions and
sometimes as solid zonal spreads.

The theory is worked out also for the Pauli-Dirac operator (anomalous
zones) and it is tested against several effects such as Aharanov-Bohm 
effect and Lamb shift. Regarding the Lamb-shift, it should be 
emphasized that the theory developed so far concerns free particles,
meaning that the Hamilton operator does not contain potential functions
$V$ attributed to other sources such as nucleus. In case of Lamb-shift,
regarding the non-existence of doublets in the spectrum of an electron
in a hydrogen atom, there is such $V$ considered which is the Coulomb
potential of the nucleus. 

The problem arising about this potential is that the
zones are not invariant with respect to multiplication with this radial
function. Therefore, this operator does not fit the zonal theory at all.
This operator, by expressing non-local interactions, 
contradicts also relativity. Einstein was able to eliminate the similar
problem arising about Newton's gravitation law by 
the theory of general relativity. Despite this intimidating 
analogy this problem will be considered in two different ways.

The first one defines the action
of the zonal Coulomb potential operator, $V^{(a)}$, 
on a zonal function 
$f^{(a)}$
by projecting 
$Vf^{(a)}$ back to $\mathcal H^{(a)}$. It turns out that
also this zonal operator is an integral operator whose smooth kernel, 
$V^{(a)}(x,y)$, spreads out (non-trivially) onto the whole space.
Thus the zonal Coulomb fields express local interactions. 
The theory will be tested against the Lamb shift
in a hydrogen atom where the nucleus has zonal Coulomb potential. 
There are no doublets in the spectrum in this case. 

In the second way the zonal theory
is extended onto curved manifolds. Thus the potential $V$ is inbuilt 
into a curved metric, like in Einstein's general relativity.

The following review meticulously goes back through the details
of the zonal theory which is just sketchily described above. 
Most of the 
mathematical formulas will be only stated with no proofs. 
The rather long mathematical details are in \cite{sz5, sz6}.

\section{Grand review of zonal theory}

{\bf (A) Brief history of renormalization.}

The confusing infinities represented by 
divergent integrals are still present 
since the very beginning of quantum electrodynamics.
The difficulties originate from concepts such as point mass
and point charge of classical electron
theory, which provided the first warning that a point 
electron will have infinite electromagnetic
self-mass; the mass $e^2/6\pi ac^2$ for a surface distribution
of charge with radius $a$ blows up for $a\to 0$. 
In his classic paper V. F. Weisskopf (1939) explains the 
trouble caused by this divergence as follows:
`` Quantum kinematics shows that the
radius of the electron must be assumed to be zero. It is easily proved
that the product of the charge densities at two different points,
$\rho (\bold r-\xi /2)\times\rho(\bold r+\xi /2)$,
is a delta-function $e^2\delta (\xi)$. In other words: if one 
electron alone is present, the probability of finding a charge
density simultaneously at two different points is zero for every
finite distance between the points. Thus the energy of the
electrostatic field is infinite as 
$W_{st}=(1/8\pi )\int(H^2+E^2)d\bold r=\lim_{a\to 0}e^2/a$.'' 

This infinity had given rise to the greatest confusion at the earliest
time of field quantization theory (cf. the pioneering work 
of Heisenberg and Pauli, published in (1929)), 
where  the above integral represents
the Hamiltonian, $H_{field}=H_f$, of the Coulomb field. 
Then this crisis further escalated without finding any 
technique by which the incessantly growing number of infinities in
the theory could have been handled. 
(A fascinating account on the development of quantum electrodinamics
from its beginning up-to the 50's can be found in {\it Silvan S. 
Schweber; QED and the men who made it: 
Dyson, Feynman, Schwinger, and Tomonaga; Princeton
University Press (1994)}. The most important original papers 
concerning this topic are 
reprinted in {\it ``Selected Papers on Quantum Electrodynamics", 
edited by Julian Schwinger, Dover Publications, 1958}).

The breakthrough on the problem
grew out of discussions at the Theoretical Physics
Conference on Shelter Island, June 2 to 4; 1947, devoted for finding a
satisfactory theoretical explanation for the so called {\it Lamb shift} 
measured between the two second level, 
$2s$ and $2p$ (resp. $2s_{1/2}$ and $2p_{1/2}$ in the Dirac theory), 
of the hydrogen atom. According to the
Schr\"odinger (resp. Dirac) equation, this second level is degenerated
and the two levels occur at the same energy. Yet, Lamb and
Retherford (1947) found that there is indeed a small separation.

The infinities experienced in QED confused all earlier 
attempts to calculate this difference. Then; at the Shelter Conference;
Kramers, Schwinger and Weisskopf, and Oppenheimer had suggested that
the possible explanation might be the shift of the energy levels 
caused by the interaction
of the electron with the radiation field. Though also this shift
comes out infinite in all theories, 
the most strongly divergent term can be identified
with an electromagnetic mass which must exist for a bound as well as for
a free electron. Therefore, this mass is already included in the 
observed mass of the electron and 
one must subtract from the theoretical expression the corresponding
expression for a free electron of the same average kinetic energy.
In this computation the desired finite quantity is produced by
the difference of two infinities.

For a clear exposition we need to explain this computation
technique in a more mathematical way. 
In the simple situation when the infinity is due to the 
infinity of the trace of
a Green kernel, $K(H_0)$, derived from the Hamiltonian $H_0$ of
a free particle, the divergent term is usually removed by perturbation.
In this scheme the particle is imagined to be placed 
in a field described, for instance,
by a potential $V$. Therefore, the Hamiltonian of the bounded particle
is $H=H_0+V$. Usually the trace of $K(H)$ is
still infinite, however, for appropriate $V$'s the difference 
$Tr(K(H)-K(H_0))$ is a finite quantity. By this technique one 
can produce, for instance, from a non-trace class
heat kernel, $e^{-tH_0}$, a relative heat kernel, 
$e^{-tH}-e^{-tH_0}$, of the trace class. This latter kernel allows
to introduce important object such as relative zeta- and eta-functions,
which could have not been defined regarding the original 
kernel \cite{mu}.

As it is well known, well defined zeta- and eta-functions describe the 
spectra of the Hamilton operators. What do then the relative zeta- and 
eta-functions describe? The answer to this question is this:
By adding $V$ to the Hamiltonian, the free particle becomes
a bounded one. This change causes shifts in the spectral lines. 
The relative functions describe, the shifts exhibited between the 
two spectra. For instance, the multiplicity
of each spectrum-element of a free Zeeman electron is infinity,
while this multiplicity drops to 2 for the bounded electron in a 
hydrogen atom. 
The relative objects describe these drastic changes (splitting) 
of the spectral lines. 

The Lamb-Retherford 
experiment demonstrated that even the doublets do not 
exist in the spectrum of hydrogen atom. The earlier computations 
retained only the Coulomb potential 
which is the main term in the interaction. To explain the Lamb-shift,
i. e. the non-degeneracy of the energy levels in the hydrogen atom, one
had to consider additional interactions between the
electron, the proton and the quantized electromagnetic field. 
The additional terms in the total Hamiltonian, which exhibit the 
energies due to these interactions, are called 
{\it radiative correction}.

Dropping of the divergent term from a mathematical 
expression is called {\it renormalization}, or, {\it regularization}. 
The idea of renormalization gained ground 
also in mathematics, used for spectral 
investigations on non-compact manifolds. This type of investigations 
started out with \cite{ops}, which was apparently motivated by 
\cite{a} and \cite{po}, written in physics.     

Regarding the Lamb shift the first non-relativistic computations using
renormalization were completed by Bethe by the end of the
Shelton Island Conference.
Then this idea was further developed by 
Dyson, Feynman, Schwinger, Tomonaga,
and others and the calculations came out in excellent agreement with
the observed value measured in the Lamb-Retherford experiment. 
Yet, this scheme for getting rid of infinities is thought to be
a theory which is imperfect on several counts. First,
it is only perturbative, second, infinities occur even if they can be
isolated and hidden \cite{v}.

In this article the {\it problem of infinities} is 
{\it approached from a completely
different angle}. Here a natural spectral decomposition of
the quantum Hilbert space will be considered and the desired finite
quantities will appear on these sectors (zones). 
This idea has fundamental effect also on the 
mathematics describing the geometry of sectors.
For instance, one can not think of points on a sector. Indeed, they
are substituted by a point-spread concept, defined by the wave functions
$\delta^{(s)}_X(Y)$ which are introduced by projecting the Dirac 
$\delta_X(Y)$-functionals onto the sector $s$. The elimination of
the point concept, which has been the fang of quantum theory since the 
very beginning, is a key idea in this
article. It terminates the infinities on the
sector-level. It terminates even Weisskopf's above described  argument,
since his $\delta_\xi$ becomes $\delta^{(s)}_\xi$, which describes 
indeed a positive
probability for the zonal particle being at two different points.
In building up this new theory one should rethink
many of the fundamental problems arising in the theory of finite
many particles. 

{\bf (B) Zeeman and Pauli operators.}

{\bf (B1) Classical Hamiltonians.}
The Hamilton function, $H_F$, and the corresponding quantum
operator, $H_Q$, of a charged particle in an electromagnetic field given
by the vector potential $\bold a$ and scalar potential $\phi$ are:
\begin{eqnarray}
H_F={1\over 2\mu}|\bold p-{e\over c}\bold a|^2+V+e\phi  ,\\
H_Q=
-{\hbar^2\over 2\mu}\Delta 
-{\hbar e\over 2\mu c\bold i}
(\bold a\cdot\nabla+\nabla\cdot\bold a)
+{e^2\over 2\mu c^2}\bold a^2+e\phi+V ,
\label{z_ham}
\end{eqnarray}
where $e$ is the charge, $\mu$ is the mass, and V is the potential
energy originated from other sources such as nucleus. 

One of the most important
version of this operator is the Zeeman operator representing a charged 
particle orbiting about the origin 
of the $(x,y)$-plane in a constant magnetic field 
$\bold B=B\partial_z$ perpendicular to the plane. 
This field is
described by the vector potential $\bold a_{(x,y)}=(B/2)(-y,x)$. There 
is also assumed that both $V$ and $\phi$ vanish, meaning 
that the particle is free and only the constant magnetic 
field is present. 

The latter assumption can be interpreted such that the free particle 
is considered in that unique inertial system where 
$\bold E=0$ and $\bold B\not =0$ holds.
Note that in the other inertial systems
the constant electromagnetic field appears with a non-vanishing 
constant electric field $\bold E\not =0$. 
Because of this fixed inertial system, all the considerations 
developed below can be linked to the relativistic notion.

It should be emphasized again that
the Zeeman zone decomposition is established
for free particles. Bounded particles in potential fields 
$V$ will be considered later, after analyzing the transitions and 
shifts the potential fields produce among the zones and spectral lines.
Also the Lamb-shift can be considered only after this analysis. 

So-far the constant magnetic field is externally added to the particle.
The hypothesis of orbiting spin in the classical theory means that the
Hamilton operator of the charged particle is supposed to be the same
as of the particle orbiting in a constant magnetic field. This
hypothesis, so to speak, internally adds a little magnet to the 
particle. The orbiting spin is exhibited by the angular momentum
operator discussed below.

An other important remark is
that the Zeeman operator will be considered 
only on the $(x,y)$-plane, thus it has the form
\begin{eqnarray}
H_Z=
-{\hbar^2\over 2\mu}\Delta 
-{\hbar eB\over 2\mu c\bold i}D_z\bullet
+{e^2B^2\over 8\mu c^2}(x^2+y^2) ,
\label{H_Z}
\end{eqnarray}
where 
$D_z\bullet =x\partial_y-y\partial_x$. In order to point out the finer 
differences between the classical theory and the one presented here,
we should make further remarks about these Hamiltonians.

{\bf (B2) Operator $H_{Zf}$.} The Zeeman Hamiltonian does not 
involve the Hamiltonian $H_f$ 
of the constant magnetic field, discussed
above. In classical field theory this term is introduced by 
renormalization. In the point-particle theory this field energy is
meaningless anyway. In the point-spread theory, however, the field
energy of the ``little magnet inside of the point-spread'' 
does have a real meaning, providing an additional term to 
the Zeeman Hamiltonian.
This combined particle-field Hamiltonian will be denoted by 
$H_{Zf}=H_Z+H_f$.

Operator $H_Q$ leaves on the 3-dimensional space, while $H_Z$
is defined over $\R^2$. The 3- and 2-dimensional versions
differ from each other just by the operator
$
-{(\hbar^2/2\mu)}\partial^2_z
$,
thus the spectral computations with respect to these two cases
can be easily compared. Also the
2-dimensional version (3) is well known. It has been intensely 
investigated, since its 
first appearance in \cite{ac}, in 
connection with the {\it Aharonov-Bohm (AB) effect} \cite{ab}. This
AB-phenomena got a lot of attention in the near past and will
be considered also in this article. A brief account
on this problem is as follows.

{\bf (B3) The AB effect} produces
relative phase shift between two electron beams enclosing a magnetic
flux even if they do not touch the magnetic field (cf. the thought
experiment performed with the AB-solenoid in \cite{ab}). This predicted
effect has no explanation in the classical mechanics and it
contradicts even relativity. Indeed, according to this theory, 
all fields must interact only locally and since the electrons do not 
reach the regions where the fields are, the effect can not be the 
production of the fields themselves. Yet, this effect was 
clearly demonstrated by the Tonomura experiments 
\cite{ton1}, \cite{ton2}. 

In their paper Aharonov and Bohm explained the predicted effect
by the ``significance of electromagnetic potentials in the quantum
theory". In the classical electrodinamics these potentials, which are
unique only upto gauge transformations, are
just convenient mathematical tools which do not have any physical
meaning. The only physical objects are the fields by which the
fundamental equations of motions can always be expressed. Nevertheless,
these potentials are needed to obtain a classical canonical formalism.

In quantum mechanics, however, the equations of motion of a particle
are replaced by the Schr\"odinger equation, which is obtained from a
canonical formalism. Thus this equation does involve the potentials,
resulting that ``in quantum theory an electron (for example) can be
influenced by the potentials even though all the field regions are
excluded from it". Then it is pointed out, in \cite{ab},
that this effect depend only on the gauge-invariant
quantity $\oint\bold a\cdot d\bold x=\int\bold H\cdot d\bold s$,
so that in reality they can be expressed in terms of the fields
inside the circuit. Yet, because of the relativistic notions 
the effect can not be interpreted as due to the fields themselves. 
There is no other way out
from this controversy but to retain a local theory by regarding the 
potential $\bold a$ as a physical variable. ``This means that we must be
able to define the physical difference between two quantum states
which differ only by gauge transformation. It will be shown in a future
paper that in a system containing undefined number of charged particles
(i. e., a superposition of states of different total charge),
a new Hermitian operator, essentially an angle variable, 
can be introduced which may give a meaning to the gauge."

In the announced second article (1964) the AB-effect is linked to 
superconductivity and super-fluidity, where the above angle variable 
was used to explain these super-phenomenas. 
Yet, a convincing mathematical model explaining the AB effect 
is still missing from the literature sofar.

The zonal theory offers a quite different
explanation to this problem. Indeed, a zonal electron is 
considered as an electron-spread, reaching every region in finite 
distance. Actually, exactly the zones contribute the desired physical 
meaning to the potential. 
The strong attachment of the zones to the potential $\bold a$ is 
revealed, for instance, by the fact that the zones can be defined 
by the quantization of the magnetic dipole moment.
Also note that the zones are not 
invariant with respect to the gauge transformations. 
These facts clearly demonstrate that
the potential reveals itself by the zones and the point-spreads defined
by the zones. 

Beyond this new explanation of the AB effect, the point-spread concept 
offers explanations also for number of other problems. Its most 
important impact is certainly the cancellation of infinities
from the theory of finite many particles.

{\bf (B4) Quantum numbers.} Term involving $D_z\bullet$ in  
operator $H_Q$ is called angular momentum operator. In what
follows, it is denoted by $\bold L$. It represents the orbiting spin 
of the particle. 
Note that the rest part of the Hamilton operator is a harmonic 
oscillator operator, $\bold O$, which commutes with $\bold L$.

By the standard method, the eigenfunctions of $\bold O$ are sought
in the form $f(r^2)\varphi_l$, where the $f$ is a radial function and
$\varphi_l$ is a homogeneous harmonic polynomial of order $l$.
For any fixed $\varphi_l$ the $\bold O$ induces a differential operator
acting on $f$. This operator depends just on $l$ and the eigenfunctions
of this radial operator appear in the form 
$f_k(r^2)=u_k(r^2)e^{-\lambda r^2/2}$,
where the $u_k$ is a $k^{th}$-order Laguerre polynomial and the 
parameter $\lambda$ is described below. Thus the multiplicity of the
eigenvalue for fixed $k$ and $l$ is the dimension of the space 
of the $l^th$-order spherical harmonics 
$H^{(l)}$. The integers $k+l$ resp. $l$
are called {\it principal resp. azimuthal quantum numbers}.

Now let also the angular momentum operator be acting. Since it commutes
with the Euclidean Laplacian, the space $H^{(l)}$ is invariant
under its action. It is well known that the $\bold iD_z\bullet$
has the eigenvalues $\{-l,-l+1,\dots ,l-1,l\}$, which are called
{\it magnetic quantum numbers}, denoted by $m$. Thus adding $\bold L$ to
the $\bold O$ causes shifts in the energies emerging on the
same level with respect to the operator $\bold O$. This shift is
called Zeeman effect whose existence was demonstrated by the
Stern-Gerlach experiment. This experiment approved the hypothesis
of magnetic dipole moment (orbiting spin, little magnet) existing 
in charged particles such as electrons.

{\bf (B5) Pauli operator.} 
Although the Schr\"odinger wave equation gives excellent  
agreement with experiment in predicting the spectral lines, yet small
discrepancies have been found which can be removed by assuming,
besides the orbiting spin, also the existence of an intrinsic angular
momentum, $\bold S$, of magnitude $-e\hbar /2\mu c$. It is well known
the many efforts made to establish an adequate spin
concept \cite{to}. The ultimate solution of this problem is due to 
Pauli, in the non-relativistic case, 
and to Dirac in the relativistic case. In
Pauli's theory the intrinsic angular momentum operator acts on
spinor fields having two components, while Dirac's theory applies
to 4-component spinor fields. In Pauli's theory the constant magnetic
field $B\partial_z$ makes the following contribution to the Hamilton
operator:
\begin{eqnarray}
   H_S=-{e\hbar\over 2\mu c}\begin{pmatrix}
     B & 0\\
     0 & -B
     \end{pmatrix} . 
\end{eqnarray}
The Hamilton operator $H_P=H_Z+H_S$, where the $H_Z$ acts on the 
2-component spinor fields as a scalar operator, 
is called {\it Pauli Hamiltonian}. A particle corresponding to $H_Z$
resp. $H_P$ is called {\it Zeeman- resp. Pauli-particle}.

{\bf (C) Mathematical modeling; Zeeman manifolds}

{\bf Zeeman operator identified with a Riemannian Laplacian.}
Most surprising feature of the Zeeman operator $H_Z$ 
is that it can be pinned down
as the Laplacian on a Riemannian manifold. As far as the
author knows, this interpretation has not been recognized 
in the literature so far. We use this Riemannian manifold as the
fundamental mathematical model describing the space-time structure on 
the quantum level. Although the metric is positive definite, 
this model fulfills the relativistic criteria.
By the first version of these manifolds single Zeeman-, or, 
Pauli-particles are modeled.

This fundamental Zeeman manifold is a Riemannian 
circle bundle over $\R^2$, defined by factorizing the center
of the 3-dimensional Heisenberg group endowed with a left-invariant 
metric. 

The Lie algebra $\bold n=\R^2\times\R=\R^3$ 
(where $\R$ is the center) of the 3D-Heisenberg
group can be described in terms of the
natural complex structure $J$, acting on $\R^2$, and the natural
inner product $\langle ,\rangle$, defined on $\bold n=\R^3$, 
by the formula
$\langle [X,Y],Z\rangle =\langle tJ(X),Y\rangle$, where the 3-vectors 
$X,Y$ and $Z$ are
in $\R^2$ and $\R$ respectively, furthermore, $t$ is
the coordinate of $Z$ in $\R$.
The map $Z\to tJ=J_Z$ associates skew endomorphisms acting on $\R^2$
to the elements, $Z$, of the center. They satisfy
the relation 
$J^2_Z=-|Z|^2id$. 
Thus the metric Lie algebra is completely
determined by the system 
\begin{eqnarray}
\{\bold n=\bold v\oplus\bold z,\langle ,\rangle ,J_Z\},
\label{nlie}
\end{eqnarray}
where $\bold v=\R^2$ and $\bold z$ are called X- and
Z-space respectively. With higher dimensional X- and Z-spaces this
system defines the Heisenberg type Lie algebras introduced by
Kaplan ~\cite{ka}. If the Clifford condition    
$J^2_Z=-|Z|^2id$
is dropped for the skew endomorphisms, the above system defines a
most general 2-step nilpotent Lie algebra. The considerations
will be extended to these general cases, however, the discussion
proceeds with the fundamental 3-dimensional case. 

Note that there
are two options, 
$J$ or $J^\prime =-J$, 
for choosing a complex structure on $\R^2$.
The two Lie algebras,
$\bold n$ and $\bold n^\prime$ are isometrically isomorphic by the
map $(X,Z)\to (X,Z^\prime =-Z)$.
 
The Lie group defined by this Lie algebra is denoted by $N$, 
furthermore, the left-invariant extension of the inner product
$\langle ,\rangle$, defined on the tangent space $T_{(0,0)}(N)=\bold n$
at the origin, is denoted by $g$. The exponential map is one-to-one 
whose inverse identifies
the group $N$ with its Lie algebra $\bold n$. Thus also the group  
lives on the same linear space $(X,Z)$. Then the
group multiplication is given by:
\begin{eqnarray}
(X,Z)(X^* ,Z^* )
=(X+X^* ,Z+Z^* +{1\over 2}
[X,X^* ]).
\label{nilprod}
\end{eqnarray}

On the linear coordinate systems 
$\big\{x^1,x^2,t\big\}$, defined by the natural basis
$\big\{E_1, E_2, e_t\big\}$, the left-invariant
extensions of the vectors $E_i;e_t$ are of the form
\begin{eqnarray}
\bold
X_i=\partial_i + \frac {1} {2} 
\langle [X,E_i],e_t \rangle  \partial_t
=\partial_i + \frac {1} {2}  \langle
J\big(X\big),E_i\rangle
\partial_t\quad ; \quad
\bold Z=\partial_t,
\label{invect}
\end{eqnarray}
where $\partial_i=\partial /\partial x^i$, $\partial_t=
\partial/\partial t$.
Then for the Laplacian, $\Delta$, acting on functions we have: 
\begin{eqnarray}
\Delta=\Delta_X+ (1+\frac 1 {4}|X|^2)
 \partial_{tt}^2+\partial_t D \bullet,
\label{debullet}
\end{eqnarray}
where $\Delta_X$ is the Euclidean Laplacian on the X-space and
$D\bullet$ means differentiation (directional derivative)
with respect to the vector field
\begin{eqnarray}
D : X \to J\big (X\big )
\label{D}
\end{eqnarray}
tangent to the X-space. 

The above Laplacian is not the desired Zeeman operator yet. 
This surprising interpretation
can be established on center-periodic Heisenberg 
groups defined by an L-periodic lattice 
$\Gamma_Z=\{Z_{\gamma L}=\gamma L|\gamma\in\mathbb{Z}\}$ 
on the center.  
Since the $\Gamma_Z$ is a discrete subgroup of isometries, 
one can consider the
factor manifold $\Gamma_Z \backslash N$ with the factor metric. The
factor manifold is a principal circle bundle over the base space 
$\bold v$ such that the circles $C_X=\pi^{-1}(X)$ over the points 
$X\in \bold v$
are of constant length $L$. Then the projection 
$\pi : \Gamma_Z \backslash N \to \bold v$
projects the inner product from the horizontal subspace
(defined by the orthogonal complement to the circles) 
to the Euclidean inner product $\langle ,\rangle$ on the X-space.

By using the Fourier-Weierstrass decomposition
\begin{eqnarray}
L^2(\Gamma\backslash N)=\oplus FW^{(\gamma)},
\label{L^2}
\end{eqnarray}
where $FW^{(\gamma)}$ consists of functions 
of the form
\begin{eqnarray}
\phi^{(\gamma)}(X,Z)=
\varphi (X)e^{\bold i\gamma 2\pi t/ L},
\label{phi}
\end{eqnarray}
the Laplacian can be established in the following particular form. 

By (8), the function spaces 
$FW^{\gamma}$ are invariant under 
the action of the Laplacian. More precisely we have:
\begin{eqnarray}
\Delta \phi^{(\gamma )}
=(\square_{(\lambda )}\varphi )e^{ 
\bold i\gamma 2\pi t/ L},\quad\quad\text{where}\\
\square_{(\lambda )}
=\Delta_X +{2 \bold i} D_\lambda\bullet -
{4\lambda^2}\big(1 + 
\frac 1 {4} |X|^2\big)\quad ,\quad \lambda ={\pi\gamma\over L},
\label{lapl}
\end{eqnarray}

and $D_\lambda\bullet =\lambda D\bullet$ means directional derivative
along the the vector field $X\to\lambda J(X)=J_\lambda (X)$. 
If $\lambda <\, 0$,
the $J$ and $\lambda$ are exchanged for $-J$ and $-\lambda$ 
respectively. Thus one can assume that $\lambda >0$.

Apart from the constant term $-4\lambda^2$, operator
$\square_{(\lambda )}$ is nothing but the Zeeman operator $H_Z$ 
described in (3). 
The surplus constant term will be identified by the field-energy of the 
constant magnetic field in the charge spread, thus the above
operator is, actually, $H_{Zf}$ described earlier. The precise  
description of identification of the Zeeman operator (3) with the
Laplacian $\square_\lambda$ acting on the invariant subspace 
$FW^{\gamma}$ is as follows. The macroscopic Zeeman operator is defined
by $\hbar =\mu =1$. Then 
$H_{Zf}=-(1/2)\square_{(\lambda )}$, 
where $\lambda =-eB/2c$. Note that particles with negative charge 
correspond to the cases $\gamma >0$, i. e., they are attached to $J$, 
while particles with positive charges are attached to $-J$.

On the quantum (microscopic) level, the periodicity $L$
and the parameter $\lambda$ are exchanged for $L_\hbar =\hbar L$
and $\lambda_\hbar =\lambda /\hbar$ respectively. This process means
nothing but scaling of the periodicity by $\hbar$. Then we have
$H_{Zf}=-(\hbar^2/2\mu )\square_{(\lambda_\hbar )}$. By scaling also the
Euclidean metric on the X-space by $\hbar/\sqrt\mu$, one has 
$H_{Zf}=-(1/2)\square_{(\lambda_\hbar )}$. In the following we proceed
with the macroscopic operator, however, the previous formulas allow
an easy transfer from the macroscopic level to the microscopic one.

By a straightforward generalization, described later, these operators 
can be introduced on higher dimensional Heisenberg groups defined
by a complex structure, $J$, on an even dimensional Euclidean space
$\R^k$.
Let $(z_1,\dots ,z_{k/2})$, where 
$z_i=q_i+\bold ip_i$ and
$\partial_{p_i}=J(\partial_{q_i})$, 
be a complex coordinate system regarding $J$. Then this system 
identifies $\R^k$ with $\mathbb{C}^{k/2}$. The circle bundle,
defined by factorization of the center $\R$, determines 
quantum operators depending just on one parameter $\lambda$. 
In the most general situation the center of a two-step
nilpotent group is factorized by a lattice, resulting a torus bundle
over the X-space. Then one arrives to operators depending on 
different parameters $\lambda_i>0$ which are defined for 
each complex coordinate plane $z_i$. As it will be explained
later, the above operators correspond to the Hamiltonians of a system 
where $k/2$ number of charged particles are circulating in a constant
magnetic field. When there is only one $\lambda$ involved, the particles
are considered to be identical upto the sign of the charge.
 
{\bf (D) Introducing the zones}

{\bf (D1) Quantum Hilbert space $\mathcal H$.} 
The quantum operator is acting on smooth functions 
primarily. This action extends to 
the $L^2$-Hilbert space $\mathcal H=L^2_{\mathbb{C}}(\R^k)$ 
of complex valued functions endowed with the inner product 
$\langle f,g\rangle=\int f\overline{g}dX$ by the Friedrics extension.
This Hilbert space can be identified with the weighted Hilbert space
$L^2_{\mathbb{C}\eta_{\lambda i}}(\mathbb{C}^{k/2})$, endowed with 
the Gaussian density $\eta_{\lambda i}$, in the following way. 

First the density is defined.
For $\lambda_i =1, \forall i,$ it has the simple form   
$
\eta =e^{-\sum z_i\overline{z}_i}
$,
while in general it is 
$
\eta_{\lambda i}=e^{-\sum\lambda_i z^i\overline{z}^i}.
$
The Hilbert space $L_{\mathbb{C}\eta_{\lambda i}}^2(\mathbb{C}^{k/2})$ 
is spanned by the polynomials written in terms of
holomorphic and antiholomorphic coordinates. The above mentioned
identification of this Hilbert space with the standard Hilbert space 
$L^2_{\mathbb{C}}$ endowed
with the standard Euclidean density $\eta =1$ is established by the map
\begin{eqnarray}
L^2_{\mathbb{C}\eta_{\lambda i}} 
\to L^2_{\mathbb{C}}\quad ,\quad
\psi\to\psi e^{-{1\over 2}\sum\lambda_i z_i\overline{z}_i}.
\label{L^2toL^2}
\end{eqnarray}

{\bf (D2) Zones defined by Heisenberg group representations.} 
The Zeeman zone decomposition of the Hilbert space 
$\mathcal H=L^2_{\mathbb{C}\eta_{\lambda i}}$ can be introduced in two 
entirely different ways. In the first way these subspaces are defined
by the semi-irreducible invariant subspaces of the 
natural reducible complex
Heisenberg algebra representation. The other definition explores 
that these subspaces are invariant under the action of the quantum
operator. The zones are introduced
by a spectral decomposition such that the spectrum of the
Hamilton operator is explicitly computed and, then, the eigen-functions 
are sorted into different classes corresponding to the distinct magnetic
states defined by the quantization of the magnetic dipole moment.

The Lie
bracket on the complex Heisenberg algebra is the restriction of the 
Poisson bracket to the set
$\{z_i,\overline{z}_i,c\}$ of   
holomorphic resp. antiholomorphic
coordinates and constant $c$. 
The real Heisenberg algebra
is hidden in this complex algebra which can be 
recovered by certain formulas established in \cite{sz4}.
The {\it representation of} this {\it complex Heisenberg algebra} is 
introduced by
\begin{eqnarray}
\rho_{\mathbb{C}} (z_i)(\psi )=(-\partial_{\overline{z}_i}+
\lambda_i z_i\cdot )\psi\quad ,\quad
\rho_{\mathbb{C}} (\overline{z}_i)(\psi )=\partial_{z_i}\psi . 
\label{heisrep}
\end{eqnarray}
This representation satisfies the well known Heisenberg relations,
however, it is not unitary on the
whole complex algebra. It becomes unitary by restricting it
to the real sub-algebra. 

The most remarkable feature of this complex representation  
is that the $\rho_{\mathbb{C}}$ is a
reducible representation on the whole Hilbert space
$\mathcal H=L_{\mathbb{C}\eta_\gamma}^2$. 
Indeed, the holomorphic subspace 
$\mathcal H^{(0)}$
spanned by the holomorphic polynomials 
$z_1^{a_1}\dots z_k^{a_k}$ 
is obviously an irreducible invariant subspace.
In the literature usually this irreducible representation is called 
complex (Fock) representation and no
thorough investigation of the whole reducible representation has been
implemented yet. 

By the first definition, the zones are introduced 
by the semi-irreducible invariant subspaces 
of this reducible representation.
This decomposition can be established by the
standard {\it Gram-Schmidt orthogonalization} such that 
the Gram-Schmidt process is applied to the series 
$G^{(a)},\, a=0,1,\dots$,
of subspaces, where $G^{(a)}$ is spanned by the
subspaces 
$\overline{z}_1^{a_1}\dots \overline{z}_k^{a_k}\mathcal H^{(0)}$ 
satisfying $a_1+\dots +a_k=a$. Clearly 
$\mathcal H^{(0)}=G^{(0)}$ holds. The next subspace,
$\mathcal H^{(1)}$, is defined as the orthogonal complement of 
$\mathcal H^{(0)}$ in
$G^{(0)}\oplus G^{(1)}$, e. t. c.,
the higher order subspaces are defined inductively.
 
Representation $\rho$ is irreducible on 
$\mathcal H^{(0)}$ and it is irreducible on the higher
order zones if and only if $k=dim(\bold v)=2$. 
In the higher dimensional cases the zones of higher order are reducible.
For instance,
the zone-functions defined by a fixed $\overline{z}_i$ form an 
irreducible invariant subspace in the higher order zones. This is an 
explanation for the meaning of the name {\it semi-irreducible}.

It should be noted that this reducibility in the higher dimensions has
a much deeper meaning. In fact, the complex dimension $k/2$ is 
interpreted as number of particles and by considering subspaces defined
by eigenfunctions satisfying certain symmetry properties with respect
to coordinate exchanges, one can define 
{\it Bosonic or Fermionic particles}
and one can introduce important concepts such as 
{\it exclusion principles}. More about this important topic can be 
found later.

{\bf (D3) Zones established by spectral decomposition 
(polarization).}  The other technique by which the 
Zeeman zones are established is based on explicit computation 
of the spectrum and eigen functions.
Next only systems involving one
parameter $\lambda$ are considered. 
If $\lambda =0$, the spectrum is the continuous spectrum of the
Euclidean Laplacian $\Delta_X$, thus we suppose $\lambda \not =0$.
In this case the spectrum is discrete whose
computation is carried out by the following ideas.

First note that the Hamilton operator $H_Z$ is the sum of the 
{\it harmonic oscillator operator}
$\bold O=\Delta_X-\pi^2\lambda^2|X|^2$ and the quantum {\it angular
momentum operator} $\bold L_\lambda =2\pi\bold i D_\lambda\bullet$.
These operators {\it commute with each other}. There is a well known
technique, using Hermite polynomials, by which 
the  eigen-functions and eigen-values of operator $\bold O$ are
explicitly determined. Then one generates a function 
space by letting the operators 
$\bold L^n_\lambda , n=1,2,\dots ,$ 
successively act on an eigenfunction determined in the first step. 
It still consists of eigenfunctions of $\bold O$
but, additionally, it is invariant also under the action of
$\bold L_\lambda$. The latter operator has the distinct eigenvalues 
$(2p-l)\lambda\bold i, p=0,\dots ,l$. 
In the last step the generated function space 
is decomposed according to these eigenvalues, in order to have
the eigenfunctions of the complete operator.

These eigenfunctions 
appear in the form $h^{(p,l-p)}(X)e^{-{\lambda\over 2}|X|^2}$, where
the $l^{th}$-order polynomial $h$ is determined by its leading term
$z_{i_1}\dots z_{i_p}\overline{z}_{j_1}\dots\overline{z}_{j_{l-p}}$.
Indexes $p$ resp. $\upsilon =l-p$ are nothing but
the number (degree) of holomorphic resp. antiholomorphic
coordinates involved into the leading term of the sought polynomial.
The eigensubspace spanned by these functions is denoted by
$\bold H^{(p,\upsilon)}$, 
where $p$ and $\upsilon$ are the so called
{\it holomorphic and antiholomorphic indexes}. 

In order to describe these eigenfunctions on the complex plane 
write the complex numbers $z$ in the polar form
$z=re^{\bold i\alpha}$.
Then the above constructed $h^{(p,\upsilon )}$ has a 
uniquely determined term $e^{\bold i(p-\upsilon )\alpha}$ which is 
multiplied by a radial function $f(r)$. I. e.,
$h^{(p,\upsilon )}(r,\alpha )=f(r)e^{\bold i(p-\upsilon )\alpha}$.
The functions of this 
form are called {\it polarized functions}. The method of 
polarization straightforwardly extends to the higher dimensional cases.

Note that polarized eigenfunctions appear in the form
$F(r)H^{(\tilde p,\tilde l-\tilde p)}$, 
where $H^{(\tilde p,\tilde l-\tilde p)}$
is an $\tilde l^{th}$-order homogeneous harmonic polynomial (spherical
harmonics). Thus it is in the eigen-subspace of magnetic dipole
moment operator $D\bullet$ belonging to the eigenvalue 
(magnetic quantum number) 
$m=\tilde p-(\tilde l-\tilde p)=2\tilde p-\tilde l$. In this case 
$\tilde p$ has the 
same meaning
as $p$ has in the first case, but it is counted only with respect to
the polynomial $H$. Also $\tilde l$ concerns strictly this polynomial.
These are very important differences between the quantum
numbers $p,l,m$ and $\tilde p,\tilde l, m$.
In standard spectrum-computation of 
electrons 
the eigenfunctions are sought exactly in the second form.
The precise calculations reveal that 
function $F(r)$ appears in the form
\begin{eqnarray}
F(r)=
u_{(\lambda ,n,\tilde l,m)}(\lambda r^2)
e^{-{1\over 2}\lambda r^2},
\label{F}
\end{eqnarray}
where
$u_{(\lambda ,n,\tilde l,m)}(t)$ is an appropriate Laguerre polynomial 
of the 
$n^{th}$-other. The eigenvalue belonging to such an eigenfunction is
\begin{eqnarray}
\mu_{(\lambda,n,\tilde l,m)}=-((4n+4\tilde p+k)\lambda+4k\lambda^2).
\label{mu}
\end{eqnarray}
Thus the relations between the two set of quantum numbers are: 
$l=2n+\tilde l\, ,\, p=\tilde p+n$, 
and $m=p-(l-p)=\tilde p-(\tilde l-\tilde p)$.

Let it be pointed out again that 
the quantum numbers are defined in the classical electron theory 
by the indexes
$n,\tilde p$ and $\tilde l$, derived from the spherical harmonics 
technique. If the electron is in the quantum state described
by the $\tilde l^{th}$-order homogeneous harmonic polynomials, 
then the eigenvalues
of the dipole moment operator are 
$m=2\tilde p-\tilde l=p-\tilde \upsilon=l-2\tilde\upsilon$, 
where $p=0,\dots ,l$. Number $m,\tilde l$ and $n+\tilde l$ are called 
{\it magnetic, azimuthal, and principal quantum numbers} respectively. 
By the above relations these classical quantum numbers can be expressed
by the ones derived from the first representation of eigen-functions 
and vice versa. 

The explicit form
\begin{eqnarray}
(2p+(k/2))\lambda+2k\lambda^2
\label{eigval}
\end{eqnarray}
of the spectrum of $H_{Zf}=-(1/2)\square_\lambda$ shows that the 
elements depend only on $p$ and they are
independent of $\upsilon$. If $-J_{\lambda_1\dots\lambda_s}^2$ 
has $s$ distinct 
eigenvalues $\lambda_1,\dots ,\lambda_s$, then the actual 
Laplace-eigenvalues depend on the corresponding $p_i$'s and they
are independent of the $\upsilon_i$'s. It follows that
{\it each spectrum element has infinite multiplicity}. 
The {\it spectral Zeeman zone 
decomposition} is defined such that 
a zone consists of
functions having the same antiholomorphic indexes 
$(\upsilon_1,\dots ,\upsilon_s)$ (which are called also {\it Zeeman
zone indexes} ($ZZI$)). 
On these
subspaces (zones) {\it the spectrum is
not just discrete but each eigenvalue 
has finite multiplicity}. 
One can prove that these spectrally
defined zones are 
equal to those defined earlier, i. e., they are the
semi-irreducible invariant subspaces with respect to the natural
complex Heisenberg algebra representations \cite{sz5}. 

If $k=2$, the multiplicity of eigenvalues is $1$ on each zone and
any 2 zones are isospectral. If 
$k>2$, {\it the elements of the spectra on
two distinct zones are the
same, however, the multiplicities are different.} One can see this by 
observing that the multiplicity of an eigenvalue depends
both on the holomorphic and antiholomorphic indexes (cf. Splitting
Theorem 2.3 in \cite{sz5}). 

The higher multiplicities of eigenvalues on the zones 
of higher dimensional complex spaces $\mathbb C^{k/2}$ allow to
introduce important concepts of multi-particles theory. 
In case of $k=2$ one can choose $J$ or $-J$
for the constructions. If one of them is attached to a particle
with positive charge, the other is attached to a particle with the
same amount of negative charge. 

In case of two particles, i. e. $k=4$
and $\{z_1,z_2\}$ is a complex coordinate system on the X-space,
there are two types of eigenfunctions in a zone. One of them
is symmetric while the other is skew symmetric with respect to
the coordinate-exchange 
$z_1\leftrightarrow z_2$. 
Pauli used such skew eigenfunctions to establish his exclusion 
principle. According to this idea, particles described by these
functions can not simultaneously be in the same quantum state (dislike
each other). They are called Fermions in the literature, which are 
characterized by the property that the
probability for being in the same quantum state is zero. 
The explanation for this phenomena is that the particles have charges
of the same sign, therefore, a repulsive force is acting between them.
In the symmetric case the particles have opposite charges and attract
each other. This explains that such particles can simultaneously be 
in the same quantum states. This
property is labeled by the name {\it Bosonic} in the literature. Both
the Fermionic and Bosonic sub-zones are invariant under the action
of the Hamilton operator and the multiplicity of each eigenvalue
is one on them. Thus each zone bears the possibility
to assign charges with the same or opposite sign to the particles.
Thus a zone can be endowed with the Pauli or anti-Pauli (Bosonic)
exclusion principle, by restricting it onto the corresponding
sub-zone.

The options for choosing signs for the charges get more
and more variegated by increasing the number, $k/2$, of the particles.
In these general cases the invariant sub-zones correspond to
dividing the complex coordinates into two classes, say 
$\{z_1,\dots ,z_r\}$ and $\{z_{r+1},\dots ,z_{k/2}\}$.
The idea is that both groups correspond to particles having 
same-sign-charges, however, the sign of the charge with respect to the 
first group is opposite to the sign in the second group. I. e.,
the skew-symmetric (Fermionic) property is valid
within both groups and they are symmetric regarding exterior
coordinate-exchanges between the two classes. 
Note that these sub-zones are not supposed to be invariant
with respect to the action of the Heisenberg group representation. 
This is clearly exhibited on the irreducible holomorphic zone 
and is present also in the classical theory where the
Heisenberg group is acting always irreducibly.

By summing up, the usage of reducible Heisenberg group representation
is one of the most important distinguishing feature of the theory
presented here. Actually, this reducibility is unavoidable. 
In fact, for a fixed azimuthal quantum number,
the distinct magnetic quantum numbers are assigned to distinct 
Zeeman zones. I. e., the zones correspond to the 
quantization of the magnetic dipole 
moment and particles in different zones are considered to be in
distinct independent magnetic quantum state. This observation clearly 
shows that there is no way to avoid the reducible Heisenberg algebra 
representation, in other words, this physical system can not be
completely described in the framework of an irreducible Heisenberg
group representation. By neglecting Zeeman zones one neglects
existing magnetic quantum states. Due to the Stone-Neumann theorem,
asserting that the infinite dimensional unitary representations
of the Heisenberg group are uniquely determined up to multiplications
with complex unities, quantum theory has always been preferring
irreducible representations.

{\bf (E) Fundamentals of point-spread geometry.}

{\bf (E1) Zonal point spreads.}
The point-spreads on a zone 
$\mathcal H^{(a)}$ 
are defined by the operator $\bold P^{(a)}$ projecting
the total Hilbert space $\mathcal H$ 
onto the zone. It turns out that these operators are integral operators
with a smooth Hermitian integral kernel $\delta_Z^{(a)}(W)$. 
It can be defined by restricting the Dirac $\delta $-kernel, 
defined by an orthonormal basis $\{\varphi_i\}$ of $\mathcal H$ by
$\delta_Z(W)=\sum_i\varphi_i(Z)\overline{\varphi}_i(W)$, onto the zones.
This restriction can be established by considering an orthonormal 
basis $\{\varphi_i^{(a)}\}_{i=1}^\infty$ on $\mathcal H^{(a)}$. 
Then the sought
projec\-tion-kernel is
$\delta^{(a)}_Z(W)=\sum \varphi^{(a)}_i(Z)
\overline{\varphi}^{(a)}_i(W)$. The sum $\sum_a\delta^{(a)}_Z(W)$
all of these zonal Dirac $\delta$'s is the total $\delta_Z(W)$.
However, this is just a formal definition of the sought zonal kernels.
Fortunately enough, all of them can be explicitly computed. 

First note that the kernel regarding the holomorphic zone is the well
known Bergman kernel. In case of a single parameter $\lambda$, it has 
the following well known form:
\begin{eqnarray}
\delta^{(0)}_{\lambda Z}(W)=
(\lambda^{k/2}/\pi^{k/2})e^{\lambda (Z\overline W-(1/2)(|Z|^2+|W|^2))}.
\label{0spread}
\end{eqnarray}
One of the new features of this article is that these projection
kernels are explicitly computed not just for the holomorphic zone
but for all zones. According to these computations, projection 
$\bold P_\lambda^{(a)}$ is of the form
\begin{eqnarray}
\bold P_\lambda^{(a)}(f)_{Z}
=\int {\lambda^{k/2}\over \pi^{k/2}}L^{((k/2)-1)}_a(\lambda |Z-W|^2)
e^{\lambda Z\cdot\overline{W}}f(W)e^{-\lambda |W|^2}dW,
\label{asproj}
\end{eqnarray}
where $L_a^{((k/2)-1)}(t)$ is the corresponding Laguerre polynomial
of $a^{th}$-order. (Formula 
$Z\cdot\overline{W}=\langle Z,W\rangle+\bold i\langle Z,J(W)\rangle$ 
is often used in establishing the formulas mentioned below.)

According to (20), a general projection kernel differs from the 
Bergman kernel just by a multiplicative Laguerre polynomial. 
More precisely we have: 
\begin{eqnarray}
\delta_{\lambda Z}^{(a)}(W)
= {\lambda^{k/2}\over \pi^{k/2}}L^{((k/2)-1)}_a(\lambda |Z-W|^2)
e^{\lambda (Z\cdot\overline{W}-{1\over 2}(|Z|^2+ |W|^2)}.
\label{aspread}
\end{eqnarray}
These zonal kernels can be interpreted such that 
a point particle concentrated
at a point $Z$ appears on the zone as an object which spreads around 
$Z$ as a wave-package with wave-function 
described by the above kernel explicitly. 

The wave-package interpretation of physical objects started out with the
de Broglie theory. This concept was finalized in the Schr\"odinger
equation. The mathematical formalism did not follow this development,
however, and the Schr\"odinger theory is built up on such mathematical
background which does not exclude the existence of the controversial
point objects. On the contrary, Weisskopf's above argument
points out that an electron must be considered as a point-object in
the Schr\"odinger theory as well. An other demonstration for 
the presence of
point particles in classical theory is the duality principle, 
stating that objects manifest themselves sometime as waves and sometime
as point particles. The bridge between the two visualizations is
built up in Born's probabilistic theory,
where the probability for a particle, 
attached to a wave $\xi$, can be found on a domain $D$ is measured by
$\int_D\xi\overline\xi$.

Although the points are ostracized from the zonal theory, 
the point-spreads 
still bear some reminiscence of the point-particles. For instance, they
are the most compressed wave-packages and all the 
other wave-functions in the 
zone can be expressed as a unique superposition of the point-spreads.
If $\xi$ is a zone-function, the above integral 
measures the probability that
the center of a pointspread is on the domain $D$. This 
interpretation restores, in some extend, the duality principle in the
zonal theory.

Function 
$
\delta^{(a)}_{\lambda Z}\overline{\delta}^{(a)}_{\lambda Z}
$
is called the density of the spread around $Z$. By this reason, 
function 
$
\delta^{(a)}_{\lambda Z}
$
is called spread-amplitude. Both the spread-amplitude and 
spread-density generate well defined measures on the path-space
consisting of continuous curves connecting two arbitrary points.
Both measures can be constructed by the method applied in constructing
the Wiener measure.
 
The point-spread 
concept bears some remote reminiscence of Heisenberg's 
suggestion (1938) for the existence of 
a fundamental length $L$, analogously
to $h$, such that field theory was valid only for distances larger than
$L$ and so divergent integrals would be cut off at that distance.
This idea has never became an effective theory, however. 
Other distant relatives of the poin-spread concept are the 
smeared operators, i. e. those suitably averaged over small 
regions of space-time, considered by Bohr and Rosenfeld in quantum 
field theory.

{\bf (E2) Global Wiener-Kac and Dirac-Feynman kernels.}
The kernels in the title are introduced in the
following unified form:
\begin{eqnarray}
e^{-\sigma tH_0}= d_{\sigma} (t,X,Y)=
\sum e^{-t\sigma\mu_i}\psi_i(X)\overline{\psi}_i(Y),
\label{gham}
\end{eqnarray}
where $\{\mu_i> 0\}$ is the discrete spectrum of the Zeeman operator
$H_Z$ and
functions in the orthonormal basis $\{\psi_i\}$ are  
eigenfunctions of $H_Z$. Furthermore, the parameters
$\sigma =1$ and $\sigma =\bold i$
correspond to the {\it Wiener-Kac resp.
Dirac-Feynman kernels}. Both kernels satisfy the condition
\begin{eqnarray}
\lim_{t\to +0}d_\sigma (t,X,Y)=\delta_X(Y),
\label{lim}
\end{eqnarray}
where  $\delta$ is the complex Dirac-delta kernel. 
The kernels attached to operator
$H_{Zf}=-(1/2)\square_{(\lambda )}$ are multiples of (22) by
$e^{2\sigma\lambda^2t}$.

Despite the infinite multiplicities, these infinite 
function series converge to
\begin{eqnarray}
d_{1} (t,X,Y)=\big(
{\lambda\over 2\pi sinh(\lambda t)}\big)^{k/2}e^{
-{\lambda\over 2}coth(\lambda t)
|X-Y|^2+\bold i
\langle X,J_\gamma (Y)\rangle},
\label{gwk}\\
d_{\bold i}(t,X,Y)=
\big({\lambda\over 2\pi\bold i sin(\lambda t)}\big)^{k/2}e^
{\bold i\{
{\lambda\over 2} cot(\lambda t)|X-Y|^2-\langle X,J_\gamma (Y)\rangle\}}.
\label{gfc}
\end{eqnarray}
 
There are numerous differences between these two kernels.
They provide the fundamental solutions
for the heat- and the Schr\"odinger equation respectively. 
The first one defines the well established 
{\it Wiener-Kac measure} on the continuous path-spaces while this
construction technique can not be carried over to the DF-kernel in 
order to construct the  {\it Feynman measure}. In fact, one arrives
(by this technique) to the well known
divergent integrals, meaning that the
approximating measures in the construction do not extend to a 
continuous complex measure on the path-space. It is known from the 
history of QED that Kac,
who tried to understand Feynman, was able to introduce a well
defined measure on the path-spaces only by the kernel $d_1$, which in
the Euclidean case ($\lambda =0$) is nothing but the Wiener measure
by which the Brownian motion is introduced. Later, in the Feynman-Kac
formulas, the Radon-Nikodym derivative of the more general Wiener-Kac
measure with respect to the Euclidean Wiener measure was determined.
This well defined measure is a very important tool in quantum theory
even today. Despite the controversies, also 
the Feynman measure is still a very powerful intuitive tool 
in QED. It should be also mentioned that none of 
these kernels is of the trace class, therefore, 
objects such as zeta functions,
eta functions, e. t. c. can be defined only by 
regularizations.

{\bf (E3) Zonal Wiener-Kac and Dirac-Feynman kernels.} 
Since the global flows leave the zones invariant,
the zonal WK- and DF-kernels can be defined by restricting the
global kernels onto the zones. The particular beauty of this theory
is that all these objects can explicitly be computed. Unlike the
global ones, both zonal kernels are of the trace class, having well 
defined partition functions. In case of a single parameter $\lambda$,
the zonal WK-kernels and the corresponding partition functions are 
of the form
\begin{eqnarray}
d_{1}^{(a)}(t,X,Z)=\label{zonwk}\\
\big({\lambda e^{-\lambda t}\over \pi}\big)^{k\over 2}
\cdot L^{({k\over 2}-1)}_a(\lambda |X-Z|^2)
e^{ \lambda (-{1\over 2} (|X|^2+|Z|^2)+ e^{-2\lambda t}
\langle X,Z+\bold iJ(Z)\rangle )},\nonumber\\
Z_1^{(a)}(t)=Trd_{1}^{(a)}(t)
={a+(k/2)-1\choose a} {e^{-{k\lambda t\over 2}}\over
(1-e^{-2\lambda t})^{k\over 2}}\, ,
\label{partwk}
\end{eqnarray}
while the corresponding zonal DF-kernels and partition functions are:
\begin{eqnarray}
d_{\bold i}^{(a)}(t,X,Z)=\label{zondf}\\
\big({\lambda 
e^{-\lambda t\bold i}
\over \pi}\big)^{k\over 2}
\cdot L^{({k\over 2}-1)}_a(\lambda |X-Z|^2)
e^{\lambda (-{1\over 2} (|X|^2+|Z|^2)+ 
e^{-2\lambda t\bold i}\langle X,Z+\bold iJ(Z)\rangle )}\nonumber\\
=e^{-{1\over 2} k\lambda t\bold i}
e^{\big(\lambda (e^{-2\lambda t\bold i}-1)\big)
\langle X,Z+\bold iJ(Z)\rangle}\delta^{(a)}_{\lambda}(X,Z),\nonumber\\
Z_{\bold i}^{(a)}(t)=Trd_{\bold i}^{(a)}(t)
={a+(k/2)-1\choose a} {e^{-{k\lambda t\over 2}\bold i}\over
(1-e^{-2\lambda t\bold i})^{k\over 2}}.
\label{partdf}
\end{eqnarray}
By these formulas the zonal zeta and eta functions can be introduced 
with no regularization. These formulas are 
beyond the scope of this review.

{\bf (F) Linking to the blackbody radiation; Solid zonal particles}

{\bf (F1) Blackbody radiation and specific heat
of solids.}
The above partition functions allow to link the evolution of the 
point-spreads driven by the zonal WK- resp. DF-flow to the blackbody 
radiation in equilibrium. Quantum theory grew out from this historic
problem, concerning the amount of
energy $U(\nu )d\nu$ radiated by the blackbody in the frequency 
range between $\nu$ and 
$\nu +d\nu$. In equilibrium the rate at which the walls emit this
frequency is balanced by the rate they absorb this frequency.
Experiments show that the $U(\nu )$ depends only on the temperature
$T$, and not on the material of which the walls are made. 

The old theory predicted that this radiation yields 
the Rayleigh-Jeans law: 
$U(\nu )d\nu\sim\kappa T\nu^2d\nu$
which contradicted the empirical curve. The empirical curve is 
fairly good described by the Wien law: 
$U(\nu )d\nu\sim\nu^3e^{-h\nu/\kappa T}$, where
$\kappa$ is Boltzmann's constant.

The controversy arising between theory and experiment was resolved by 
Planck by the hypothesis that the energy
attached to frequency $\nu$ is restricted to the integral multiple
of the basic unit $h\nu$, i. e. $E_n=nh\nu$, where $n$ is
any positive integer number. Furthermore, the probability that the 
wall emits-absorbs an energy-quanta $E_n$ is 
$W(n)\sim e^{-E_n/\kappa T}=e^{-nh\nu/\kappa T}$.
Thus by normalization we have:
$W(n)=e^{-nh\nu/\kappa T}(1-e^{-nh\nu/\kappa T})$.
Then a simple calculation yield that the mean energy is
$\tilde E=\sum E_nW(n)=h\nu e^{-h\nu/\kappa T}/
(1-e^{-h\nu /\kappa T})$,
from which the Planck distribution
$U(\nu )\sim h\nu^3 e^{-h\nu/\kappa T}/
(1-e^{-h\nu /\kappa T})$ follows. The latter distribution is the exact
form of the Wien law. Let it also be mentioned that the Rayleigh-Jeans 
law is derived by assuming the equipartition of the energy.

Einstein proposed to adopt Planck's idea to solids in order  
to explain the experimental curve describing the specific heat,
defined by the change in energy with temperature, of materials.
This curve approaches zero at absolute zero, and rises asymptotically
to $\kappa $ per atom at high temperatures. In contrast to the
blackbody radiation, where all possible frequencies can occur,
the materials have only one frequency, which is the characteristic 
frequency of the substance. Then
$\partial_T\tilde E=(h\nu)^2 e^{-h\nu/\kappa T}/\kappa T^2
(1-e^{-h\nu /\kappa T})^2$,
which formula clearly fits the curve of specific heat.

{\bf (F2) WK-flow linked to blackbody radiation of solids.}
In case of the WK-kernel the microscopic heat flow formula can be 
linked to the blackbody radiation of solids by the substitutions
$\nu =1, T=1/t, \kappa =2\pi\mu/\lambda$, where $\lambda =|e|B/2c$.
The possible integral values for the quantized energy are $n=2p+(k/2)$,
where the $p$ can take any integral value between $0$ and $\infty$.
For $k=2$ the multiplicity of eigenvalues on a zone is $1$ and
the probability $W_1(n)$ is defined by
$W_1(n)\sim e^{-E_n/\kappa T}=e^{-nh\lambda t/2\pi\mu}$, thus 
one gets the exact probability by dividing this by 
$
Z_1^{(a)}(ht/2\pi\mu )
$. These quantities are independent from the zone-index $a$. 
If $k>2$, then multiplicities occurs and these formulas should be
defined accordingly. These quantities depend on the zone index via
the combinatorial factor in (27). Since these general formulas 
can be easily recovered by multiplying with the combinatorial factor
and the appropriate exponentiation with $k/2$, we proceed with
the case $k=2$. 
A simple calculation yields 
\begin{eqnarray}
\sum_nE_ne^{-E_n\lambda t/2\pi\mu}
=-{2\pi\mu\over\lambda}\partial_t
Z_1^{(a)}(ht/2\pi\mu )=\label{aveever}\\
Z_1^{(a)}(ht/2\pi\mu )\big(h+2h
{e^{-2h\lambda t/2\pi\mu}\over 1-
e^{-2h\lambda t/2\pi\mu}}\big),
\nonumber
\end{eqnarray}
thus the average energy is
\begin{eqnarray}
\tilde E_1(T)=h+2h
{e^{-2h/\kappa T}\over 1-
e^{-2h/\kappa T}}.
\label{aveener2}
\end{eqnarray}
Note that the second term is the average energy in the blackbody
radiation of solids with the characteristic frequency $\nu =2$.
Thus Einstein's specific heat formula is
\begin{eqnarray}
\partial_T\tilde E_1(T)={(2h)^2 e^{-2h/\kappa T}\over\kappa T^2
(1-e^{-2h/\kappa T})^2}.
\label{specheat}
\end{eqnarray}
Also note that at $t=0$, which corresponds to $T=\infty$, the WK-flow 
starts out with the equipartition of energy, i. e. with the 
Rayleigh-Jeans
law what is turned to the Wien law for all $t>0$, or, $T<\infty$.
At the start all the quantities $Z_1, E_1, \partial_TE_1$ are 
infinities.

{\bf (F3) Stable (solid) charge spreads.}
The same computations lead to completely different interpretations
in case of the DF-kernel.
The analogously defined complex measure $W_{\bold i}$ is 
interpreted as probability amplitude for the energy emission-absorption
in the blackbody radiation in equilibrium. Then the average energy
and specific heat amplitudes are
\begin{eqnarray}
\tilde E_{\bold i}(T)=h+2h
{e^{-2h\bold i/\kappa T}\over 1-
e^{-2h\bold i/\kappa T}}\,\, ,\,\,
\partial_T\tilde E_{\bold i}(T)={(2h)^2\bold i\, e^{-2h\bold i/\kappa T}
\over\kappa T^2(1-e^{-2h\bold i/\kappa T})^2}.
\label{fspecheat}
\end{eqnarray}
In the following these formulas are analyzed in terms of $t$. 
The DF-specific heat does 
not bear the same physical meaning attributed to the WK-specific heat. 
In fact, the limit of $|\partial_T\tilde E_{\bold i}|$ is $\kappa$, 
both at $0$ and $\infty$. 
On the other hand, an interesting interpretation 
can be given for the other functions. 

Both $Z_{\bold i}(t)$ and $\tilde E_{\bold i}(t)$ 
are L-periodic functions, where $L=\kappa/\hbar$. 
At the lattice points $nL/2$,
where $n$ can be arbitrary integer, they have limit at the infinity,
furthermore, the density functions 
$|Z_{\bold i}|^2$ and $|\tilde E_{\bold i}|^2$ attain their minimum 
exactly at points $(n+(1/2))L/2$ of the periodicity intervals.
Thus starting from the left endpoint $P_{n}$ on an interval the 
infinite average energy drops to its minimum at the quarter point
$P_{n+1/4}$, then it grows up to the infinity at the midpoint 
$P_{n+1/2}$.
On the second half it drops to the minimum at the three-quarter point
$P_{n+3/4}$, then it 
increases to the infinity on the right endpoint $P_{n+1}$. 
This phenomena is exhibited also by the
function 
$
d_{\bold i}^{(a)}(t,X,X)=
{(\lambda 
e^{-\lambda t\bold i}
/ \pi)}
e^{\lambda ( 
e^{-2\lambda t\bold i}-1)|X|^2}
$
pointwise, since
$
|d_{\bold i}^{(a)}(t,X,X)|^2
$
attains maximum, for any fixed $X$, exactly at the end- and 
mid-points of the
periodicity intervals and the minimum is taken at the quarter points
$P_{n+1/4}$ and $P_{n+3/4}$. 

The tension amplitude and density at $X$ are defined by
\begin{eqnarray}
\tau^{(a)} (X)=\partial_t
d_{\bold i}^{(a)}(t,X,X)\quad
\text{and}\quad |\tau^{(a)}(X)|^2
\label{tau}
\end{eqnarray}
respectively. Thus there is no tension in
the minimum state, while it is infinity at the maximum places. 
This and the above considerations suggest
that the uniquely determined spreads
\begin{eqnarray}
d^{(a)}_{\bold i}(n+(1/4))L,X,Z)=-\bold ie^{-2\lambda X\overline Z}
\delta^{(a)}_\lambda (X,Z)\label{min1}\\
=-{\lambda \bold i\over\pi}L_a(\lambda |X-Z|^2)
e^{-\lambda ((1/2)|X+Z|^2+\bold i\langle X,J(Z)\rangle )}\, ,\nonumber\\
d^{(a)}_{\bold i}(n+(3/4))L,X,Z)=\bold ie^{-2\lambda X\overline Z}
\delta^{(a)}_\lambda (X,Z)\label{min2}\\
={\lambda \bold i\over\pi}L_a(\lambda |X-Z|^2)
e^{-\lambda ((1/2)|X+Z|^2+\bold i\langle X,J(Z)\rangle )}
\nonumber
\end{eqnarray}
at the quarter points should be considered as the stable charge spreads.
If the particle is in a higher average energy state then, due to the
tensions, it drops down to a stable minimum state. 

Also note that the Dirac
delta spreads represent the Rayleigh-Jeans law and the DW-flow
moves the spread down to the minimal state, which represents the Wien 
law. The introduction of stable zonal charge spreads seems to terminate
the de Broglie-Schr\"odinger's waves from the theory. The contradiction
can be resolved, however, by the duality principle, asserting that
charged materials behave sometimes as waves and sometimes as particles
being in the stable charge-spread state.

This interpretation provide a plain explanation not just for 
the problem why stable charge spreads do not blow up but it predicts 
also the spontaneous emission of accelerated electrons. 
It is well known that such electrons radiate
even when no light is incident. The explanation for this phenomena
is as follows. The accelerations
displace the electron spread from the stable minimal state and it moves 
back by radiating the energy excess gained by the acceleration.

{\bf (G) Zonal WK- and Feynman-measures.}

{\bf (G1) Construction of WK- and Feynman-measures.} 
The existence of zonal WK-measures on path-spaces follow  
from the well
defined global WK-measure. Nevertheless, it is rather striking that 
the zonal DF-kernels are
not just of the trace class but they well define, regarding any 
zone, a complex measure
$dw_{\bold ixy}^{T(a)}(\omega)$ 
on the space $\mathcal P_{xy}^T$ of continuous paths
$\omega :[0,T]\to\R^k$ connecting points $x$ and $y$. (This 
path-space is topologized with the topology of uniform convergence. 
By the one-point compactification 
$M=\R^k\cup \infty$ the $\mathcal P_{xy}(M)$ becomes a compact 
topological space.) Well defined measures, 
$d\nu_{xy}^{T(a)}(\omega)$ and $d(\nu\overline\nu)_{xy}^{T(a)}$, 
can be constructed also by the zonal spread amplitudes 
$\delta^{(a)}$ 
and densities 
$\delta^{(a)}\overline{\delta}^{(a)}$ respectively.

All these measures are constructed by the same ideas 
the Wiener measure, regarding the Euclidean heat kernel $E(t,p,q)$, 
was established. This construction starts out with the elementary fact 
that the Borel $\sigma$-algebra
(generated by the open sets) on $\mathcal P_{xy}^T(M)$ 
can also be generated by the fibred sets
$\rho_{\bold t}^{-1}(B)\subset\mathcal P_{xy}^T(M)$, where 
$\bold t=0<t_1<\dots <t_n<T$ is a fixed subdivision, $\rho_{\bold t}:
\mathcal P^T_{xy}(M)\to
M^n=M\times\dots\times M$ is the evaluation map defined by
$\rho_{\bold t}(x)=(x(t_1),\dots ,x(t_n))$, and $B$ is a Borel subset
of $M^n$. Measure $w_{Exy}$ on a fibred set 
$\rho_{\bold t}^{-1}(B)$ is defined by
\begin{eqnarray}
w_{Exy}\big(\rho_{\bold t}^{-1}(B)\big)=\label{wmeas}\\
\int_BE(t_1,x,m_1)E(t_2-t_1,
m_1,m_2)
\dots E(T-t_n,m_n,y)dm_1\dots dm_n.
\nonumber
\end{eqnarray}
By classical results, such as Riesz' theorem (concerning the measure
representation of bounded linear functionals on the Banach space of
continuous functions defined on a compact metrizable space where
the Banach 
norm is defined by $sup |f|$) and the Stone-Weierstrass theorem 
(asserting that
the curves $\rho_{\bold t_n}^{-1}(x_1,\dots ,x_n)$ corresponding
to $T{1\over n}<\dots <T{n-1\over n}$ are dense in 
$\mathcal P^T_{xy}(M)$), 
this construction determines a complex countably
additive regular Borel measure 
$w^T_{Exy}$ on $\mathcal P^T_{xy}(M)$, satisfying
\begin{eqnarray}
w^T_{Exy}\big(\mathcal P^T_{xy}(M)\big)=E(T,x,y).
\label{totmeas}
\end{eqnarray}

The same idea works out for any of the kernels mentioned above. The
construction is completed by proving the uniform 
boundedness of the approximating measures. This is established by a
particular integral formula in \cite{sz5}.
Although these measures are constructed on the total
space $\mathcal P_{xy}(M)$, the zonal measures are concentrated on 
curves
which can be described by functions constituting the zones. We do not
consider these technical details in this paper.

{\bf (G2) Zonal Feynman-Kac type formulas, Feynman's only 
stopwatch.}
It should also be pointed out 
that the zonal Euclidean heat- or DF-kernels
are not defined because the zones are not invariant with respect to the
action of the Euclidean Laplacian $\Delta_X$. Therefore, the original
Feynman-Kac formulas can not be formally carry over to the zones. 
This is why $d\nu_{xy}^{T(a)}(\omega)$ and 
$d(\nu\overline\nu )_{xy}^{T(a)}$ are introduced. 
In the zonal Feynman-Kac type formulas the Radon-Nikodym derivative 
of $dw_{\sigma xy}^{(a)}$ with respect 
to the latter measure are explicitly computed. 
The Radon-Nikodym derivative
of the zonal Feynman measure with respect to the zonal Wiener-Kac
measure is also established. These zonal FK-type formulas
have the following explicit form:
\begin{eqnarray}
d_{1}^{(a)}(T,x,y)=\int_{\mathcal P^T_{xy}} e^{-{kT\over 2}
-2\int_0^T|\omega (\tau )|^2d\tau }d\nu_{xy}^{T(a)}(\omega ),
\label{radon1}
\\
d_{\bold i}^{(a)}(T,x,y)=\int_{\mathcal P^T_{xy}} e^{(-{kT\over 2}
-2\int_0^T|\omega (\tau )|^2d\tau )\,\bold i}d\nu_{xy}^{T(a)}(\omega ).
\label{radon2}
\end{eqnarray}
The Radon-Nikodym derivative of $d\nu_{xy}^{T(a)}$ resp. the zonal 
Feynman measure with respect to the zonal WK-measure are: 

\begin{eqnarray}
d\nu_{xy}^{T(a)}(\omega )=e^{{kT\over 2}
+2\int_0^T|\omega (\tau )|^2d\tau }dw_{\bold 1xy}^{T(a)}
(\omega ),\label{rad3}
\\
dw_{\bold ixy}^{T(a)}(\omega )=e^{({kT\over 2}
+2\int_0^T|\omega (\tau )|^2d\tau )(1-\bold i)}dw_{\bold 1xy}^{T(a)}
(\omega ),\label{radon4}
\end{eqnarray}

which establish the most direct connection between the three measures.
The above formulas have the very same form with respect to any kernel 
$e^{-t\sigma H_Z}$, where the $\sigma$ is
arbitrary unit complex number. One should just substitute
$\sigma$ for $\bold i$, however, the last $\bold i$ in the firs
equation should be left alone, obviously. 

The zonal Feynman measure of the whole set of curves connecting 
$x$ and $y$ is
\begin{eqnarray}
w_{\bold ixy}^{T(a)}\big(\mathcal P_{xy}^T(M)\big)=
d_{\bold i}^{(a)}(T,x,y).
\label{ftot}
\end{eqnarray}
By an {\it intuitive interpretation of Feynman}, 
where he referred to the global kernel and measure,
the motion of a particle 
(electron, photon, e. t. c.) from $x$ to $y$ is timed by a ``stopwatch"
whose hand starts rapidly turning when leaving $x$ and 
it is stopped when 
the particle arrives to $y$. The ``hand" of the stopwatch is
considered as a complex unit number
and this timing is performed along each continuous curve connecting 
$x$ and $y$. Suppose that the particle moves from $x$ to $y$ during
the same outer time $T$, which is not the time measured by the 
stopwatch. Speaking hypothetically,
the particle is moving along each of 
the paths $\omega :[0,T]\to \R^k\, ,\omega (0)=x\, ,\omega (T)=y$
connecting $x$ and $y$. Thus, by the
finite arrow (hand of the stopwatch), each curve is 
represented by a unit 
complex number at $y$.
On the other hand, there is also a complex measure defined on the 
path-space, by which these final arrows are integrated, producing
a single final arrow (complex number). 
This unique complex number is
then $d_{\bold i}(T,x,y)$ which is called {\ probability amplitude}.
The positive real number $d_{\bold i}\overline{d}_{\bold i}$ defines
the {\it probability density} at $y$.

For constructing the measure on the path-space, Feynman used
the global kernel $d_{\bold i}$ in the same way how $d_1$ is used
for constructing the well defined Wiener-Kac measures. In case of the
Feynman measures, however, the approximating measures diverge and
they do not extend into a continuous complex measure defined on the
Borel sets of the path-space. That is why the above intuitive 
(yet very beautiful) idea is considered to be mathematically imperfect. 

On the zonal setting, however, this idea works out
perfectly and both the zonal
Feynman measure and the {\it zonal Feynman stopwatch} can be explicitly
determined. The {\it zonal probability density} is defined by 

\begin{eqnarray}
\rho^{T(a)}_{xy}=\pi^{k/2}
d_{\bold i}^{(a)}(T,x,y)
\overline{d}_{\bold i}^{(a)}(T,x,y),
\label{prob}
\end{eqnarray}

meaning that, for a Borel set
$B$, the integral $\int_B\rho^{(a)}_{x,y}db$ measures the probability
that the point spread about $x$ can be caught, at the time $T$, 
among the point spreads whose centers are on the set $B$. 
The probability regarding the whole space is always $1$, regardless
the time $T$. Feynman's stopwatch, the turning unit complex number,
is explicitly determined by the zonal Radon-Nikodym derivative (density)

\begin{eqnarray}
{dw_{\bold ixy}^{T(a)}\over d\nu_{xy}^{T(a)}}(\omega )=
e^{(-{kT\over 2}-2\int_0^T|\omega (\tau )|^2d\tau )\,\bold i},
\label{radnik}
\end{eqnarray}
which is the integrand in (40). The fascinating thing is that Feynman
has only one stopwatch for all of the zones. The 
zone-depending object is the measure $d\nu^{T(a)}_{xy}$ by which the
arrows at $y$ are integrated. Also the densities in the zonal
Radon-Nikodym derivatives (39), (41), and (42) are independent from the
zones. 

The very same statement can be established for {\it probabilities
defined by a normalized zonal wave function} $\psi^{(a)} (t,X)$
satisfying the Schr\"odin\-ger equation. Such a function
is the uniquely determined extension of the initial function 
$\psi^{(a)}(0,X)$. The extension is defined by the convolution formula 
$\psi^{(a)}(t,X)=d^{(a)}_{\bold i}(t,X,Z)*_Z\psi^{(a)}(0,Z)$. The
probability concerning the density function 
$\psi^{(a)}\overline{\psi}^{(a)}$  is interpreted as the
likely-hood that the zonal object described by the
wave function can be caught on a Borel set $B$ at the time $T$. 
Like the first one, also this probability satisfies 
the conservation law.
All these statements are particular exhibitions of the theorem 
asserting that {\it
the Feyn\-man-Dirac zonal flows define a unitary 
semigroup, $U_t^{(a)}$, on each zone}. Thus one has a unitary 
semigroup, $U_t=\oplus_aU_t^{(a)}$, on the whole $L^2$ function space. 

As there is 
explained earlier, the probabilistic theory is the bridge connecting
the de Broglie-Schr\"odinger waves with the particles defined
by the stable zonal charge spreads.

{\bf (H) Zeeman manifolds with higher dimensional centers.}

The mathematical model for interpreting the Zeeman operator
as the Laplacian on a Riemannian manifold has been, so-far, 
a Riemannian circle bundle, defined by factorizing the centers on 
Heisenberg groups which is endowed with a left invariant metrics.
This idea works out also on metric two-step nilpotent Lie groups, which
are rudimentary described in section (C), whose center
$\bold z$ is factorized by a lattice $\Gamma_Z$. 
This center is considered as
an abstract higher dimensional space such that an element $Z\in\bold z$
is identified with the endomorphism $J_Z:\bold v\to\bold v$ and 
its natural inner product is defined by
$\langle Z_1,Z_2\rangle =-Tr(J_{Z_1}\circ J_{Z_2})$.
Formulas (5)-(13)
in (C) apply also to these general cases, just the Laplacian (13)
appears in a slightly different form. 
Upto isomorphism, the Lie algebra of such a group is uniquely 
determined by a linear space, $J_{\bold z}$, of skew endomorphisms 
acting on the Euclidean space $\bold v$. Two 2-step nilpotent groups
are isometrically isomorphic if and only if the corresponding 
endomorphism spaces are conjugate. 

The rather large 
class of Riemannian torus bundles introduced in this way are  
called also Zeeman manifold. Below also particular 
Zeeman manifolds are introduced. It is
remarkable that for the so called Clifford-Zeeman manifolds even 
classification can be implemented, which may be used for classifying
the charged particles investigated in this theory. 
 
The Laplacian on the Riemannian group $(N_{J_{\bold z}},g)$, defined
by the endomorphism space $J_{\bold z}$, 
has the explicit form:

\begin{eqnarray}
\Delta=\Delta_X+\Delta_Z+\frac 1 {4} \sum_{\alpha,\beta =1}^r 
\langle J_\alpha
\big (X\big),J_\beta\big (X\big)\rangle
 \partial_{\alpha\beta}^2
+\sum_{\alpha =1}^r\partial_\alpha D_\alpha \bullet,
\label{delta}
\end{eqnarray}
which leaves the function spaces $FW^{(\gamma )}$ spanned by
the functions of the form 
$\Psi^{(\gamma )}(X,Z)
=\psi (X)e^{2\pi\bold i\langle \mathcal Z_\gamma ,Z\rangle}
=\psi (X)e^{2\bold i\langle Z_\gamma ,Z\rangle}
$,     
for all lattice points 
$\mathcal Z_\gamma\in\Gamma_Z$ 
(resp. $ Z_\gamma\in \pi\Gamma_Z$), 
invariant.
Its action on such a 
function space can be described in the form
$\Delta (\Psi^{(\gamma )})(X,Z)=\square_{(\gamma)}(\psi )(X)
e^{2\pi\bold i\langle\mathcal Z_\gamma ,Z\rangle}$,  
where operator $\square_{(\gamma )}$, acting on 
$L^2(\bold v)$, is of the form
\begin{eqnarray}
\square_{(\gamma )}
=\Delta_X + 2\pi\bold i D_{(\gamma )}\bullet -4\pi^2
\big(|\mathcal Z_\gamma |^2 + 
\frac 1 {4} |J_{\mathcal Z_\gamma}(X)|^2\big)\label{lapla} \\
=\Delta_X + 2\bold i D_{ Z_\gamma }\bullet -
4\big(| Z_\gamma |^2 + \frac 1 {4} |J_{Z_\gamma}(X)|^2\big).
\nonumber
\end{eqnarray}
 
Thus the Zeeman operator appears on the invariant
subspaces defined by the Fourier-Weierstrass decomposition. The 
spectral investigations on these manifolds are reduced to investigate
this operator on each Fourier-Weierstrass subspace separately. 

The particles represented
by these Riemannian torus bundles are called {\it Zeeman 
molecules}. A physical interpretation of factorization by 
the lattice $\Gamma_{Z}=\{\mathcal Z_\gamma\}$ is that there is a 
quantization considered also on the space of torque-axes $Z$ of the 
magnetic dipole moment.

There are special Z-molecules, defined by particular endomorphism 
spa\-ces, which are particularly interesting.
The {\it Heisenberg-type} or {\it Cliffordian
endomorphism spaces} are attached to Clifford modules 
(representations of Clifford algebras).
They are characterized by the property $J^2_Z=-|Z|^2id$, for all
$Z\in \bold z$, \cite{ka}. The corresponding molecules are called
{\it Clifford-Zeeman molecules}.
The well known {\it classification} of Clifford modules
provides classification also for the Clifford endomorphism
spaces and molecules. A brief account  on this 
classification theorem is as follows.

{\it If $r=dim(J_{\bold z})\not =3(mod 4)$, then 
there exist (up to equivalence) exactly one
irreducible H-type endomorphism space acting on a $\R^{n_r}$,
where the dimension $n_r$, depending on $r$, 
is described below. This endomorphism space
is denoted by $J_r^{(1)}$. If $r=3(mod 4)$, then there exist 
(up to equivalence) exactly
two non-equivalent irreducible H-type endomorphism spaces acting on
$\R^{n_r}$ which are denoted by 
$J_r^{(1,0)}$ and
$J_r^{(0,1)}$ 
respectively. They are connected by the relation 
$J_r^{(1,0)}\simeq
-J_r^{(0,1)}$.
 
The values $n_r$ corresponding to
$
r=8p,8p+1,\dots ,8p+7
$
are

\begin{eqnarray}
n_r=2^{4p}\, ,\, 2^{4p+1}\, , \, 2^{4p+2}\, , \,
2^{4p+2}\, ,
\, 2^{4p+3}\, ,\, 2^{4p+3}\, , \, 2^{4p+3}\, , \,
2^{4p+3}.
\label{cliff}
\end{eqnarray}

The reducible Cliffordian endomorphism spaces can be built up by these
irreducible ones. They are denoted by 
$J_r^{(a)}$ resp. $J_r^{(a,b)}$.
The corresponding Lie algebras are denoted by 
$\frak h^{(a)}_r$ resp. 
$\frak h^{(a,b)}_r$. In the latter case the X-space 
is defined by the $(a+b)$-times product 
$\R^{n_r}\times\dots\times\R^{n_r}$
such that on the last $b$ component the action of a $J_Z$ is defined by 
$J^{(0,1)}_Z$ 
and on the first $a$ components this action is defined by
$J^{(1,0)}_Z$. In the first case this process should be applied only on 
the corresponding $a$-times product.} 

In a Clifford endomorphism space each endomorphism anticommutes with
all perpendicular endomorphisms. In other words, all endomorphisms
are anticommutators. A more general concept can be 
introduced by  the {\it anticommutative endomorphism spaces} 
where all endomorphisms are anticommutators.
They can be built up, in a non-trivial way, by Clifford endomorphism
spaces. Roughly speaking, a CZ-molecule is the compound of
irreducible molecules of the same type while an {\it anticommutative
Z-molecule} is the compound of CZ-molecules of different types in 
general.

Originally, the metric groups   
$(N_J,g)$ were used, in many different ways,  for constructing 
isospectral Riemannian metrics with different local geometries. The
author's results regarding such constructions are published in
\cite{sz1,sz2,sz3,sz4} which contain also
detailed history about this topic. 
These examples 
include isospectral pairs of metrics on
ball$\times$torus-, sphere$\times$to\-rus-, ball-, and sphere-type
manifolds. Among these examples the most striking are
those constructed both on sphere- and sphere$\times$torus-type
manifolds. One of the metrics in the isospectral pair
is homogeneous while the other one is not even locally homogeneous.
These isospectrality constructions are implemented such that on some
of the irreducible subspaces $\R^{n_r}$ the endomorphism spaces
$J_r^{(1,0)}$ 
(resp. $J_r^{(0,1)}$) are switched to 
$J_r^{(0,1)}$ 
(resp. $J_r^{(1,0)}$). It turns out that the Riemannian space, 
resulted by this switching, has a completely different local geometry,
yet, the considered domains in the original and the
new Riemann spaces are isospectral.  
Endomorphism spaces $J_r^{(1,0)}$ and $J_r^{(0,1)}$ are considered 
to be representing irreducible CZ-particles having opposite charges. 
Thus the isospectrality theorem can be physically interpreted as 
follows: 

{\it By charging some of the irreducible CZ-particles 
in a CZ-molecule by the same amount of the opposite
charge the spectra of the considered domains 
remain the same, however, the local
geometry is drastically changed in general.}

Most of these isospectrality statements are established by constructing
intertwining operators, while some are proved by explicit
computations of the spectrum. These computations are
different from the one developed for the Zeeman zones.
They are, rather, 
the relative of techniques applied in physics for  
computing the spectra of charged
particles in a Coulomb potential field.
If such potential is present the eigenfunctions can be sought just 
in the form $F(r)H^{(p,l-p)}$, described in section (D3). There is
also explained that the Zeeman zones can not be constructed by 
these eigenfunctions.

{\bf (I) The Pauli-Dirac operators.}

{\bf (I1) Introducing the PaDi-operator.}
In this review the Pauli-Dirac operator is considered only on
the plane, $\mathbb C$, of complex numbers $z=(z_1,z_2)$. In
higher dimensions this operator is the sum of operators
defined on the complex coordinate planes.
These operators are square roots of the Hamilton operators, 
introduced by means of matrices 

\begin{eqnarray}
   \sigma_1={1\over\sqrt 2}
     \begin{pmatrix}
     0 & 1\\
     1 & 0
     \end{pmatrix} 
   +{\bold i\over\sqrt 2}
     \begin{pmatrix}
     0 & -1\\
     1 & 0
     \end{pmatrix} \,\, ,\,\, \sigma_2=\overline{\sigma}_1\,\, ,\,\,
   \sigma_0 =\begin{pmatrix}
     1 & 0\\
     0 & -1
     \end{pmatrix}.
\label{paulim}
\end{eqnarray}

Matrices $\sigma_1$ and $\sigma_2$, 
which are built up by the well known 
Pauli spin matrices, 
are called
canonically conjugate spin matrices. They satisfy
the commutativity relations

\begin{eqnarray}
\sigma_i\sigma_j+
\sigma_j\sigma_i=2\delta_{ij}
     \begin{pmatrix}
     1 & 0\\
     0 & 1
     \end{pmatrix} .
\label{clifpr}
\end{eqnarray}

Such pairs can be defined 
for any pair $\{u_1,u_2\}$ of perpendicular unit complex numbers, 
which define the coefficients before the Pauli matrices. 
The above matrices correspond to
$u_1=(1+\bold i)/\sqrt 2, u_2=\overline u_1$.
Any of such pairs $(\sigma_1,\sigma_2)$ is appropriate
to establish a PaDi-operator. Note that $u_2$ is not the conjugate of 
$u_1$ in general; this is true just for the above matrices. I. e.,
canonically conjugate spin matrices are generated by perpendicular
and not conjugate complex numbers. 

The PaDi-operator is defined by 

\begin{eqnarray}
\mathcal{PD}={1\over \sqrt 2}\sum_{j=1}^2\sigma_j
(\bold i\partial_{z_j}-\bold a^j)+
2\sigma_0\lambda =
\label{paudi}\\
=     \begin{pmatrix}
2 \lambda
      & 
{1+\bold i\over\sqrt 2}(2\partial_{\overline z}-\lambda z)
\\
{-1+\bold i\over\sqrt 2}(2\partial_{z}-\lambda\overline z)
      & 
-2 \lambda
     \end{pmatrix} .
\nonumber
\end{eqnarray}

This PaDi-operator is attached to the
Hamiltonian $H_{Zf}$.
The corresponding operator attached to $H_Z$ is defined by 
omitting the second
term, $\lambda\sigma_0$, in (51). They
are distinguished by the denotations 
$\mathcal{PD}_{Zf}$ and $\mathcal{PD}_Z$.

The PaDi-operator acts on $\mathbb C^2$-valued functions, called
2-component spinor fields, which are written
in the form $\phi=(\varphi_1,\varphi_2)$.
The inner product of spinor fields
$\phi$ and $\gamma$ is defined by
 
\begin{eqnarray}
\langle\phi,\gamma\rangle =\int_{\R^2} \sum_i\varphi_i(X)
\overline{\gamma}_i(X)dX.
\label{innprod}
\end{eqnarray}
 
The corresponding $L^2$ spinor Hilbert space is denoted by $\mathcal S$.

In order to compute the squared operators,
we express $\mathcal{PD}$ in a more explicit form. Since 
the vector potential is of the form
$\bold a=\lambda (-z_2,z_1)$, thus

\begin{eqnarray}
\mathcal{PD}_Z(\phi )=(\mathcal{D}_1(\varphi_2),
\mathcal{D}_2(\varphi_1)) \, ,\quad \text{where}\label{PD}\\
\mathcal{D}_1=
{1+\bold i\over \sqrt 2}(2\partial_{\overline z}-\lambda z)
\quad ,\quad
\mathcal{D}_2={-1+\bold i\over\sqrt 2}(2\partial_z+\lambda\overline z).
\label{paudi2}
\end{eqnarray}
Though the component operators act only on smooth functions,
by Friedrics extension, their action extends to the function 
space $L^2_{\mathbb C}$. By (14), 
this space is isomorphic to the weighted 
Hilbert space $\mathcal H=L_{\mathbb C\eta}^2$, 
defined by the Gaussian density
$
\eta_\lambda =e^{-\lambda |X|^2}.
$
This isomorphism defines an isomorphism also between $\mathcal S_\eta$
and $\mathcal S$.

On the weighted Hilbert space operators (54) appear in the form

\begin{eqnarray}
\mathcal{D}_1=
\sqrt 2(1+\bold i)
(\partial_{\overline z}-\lambda z)\quad ,\quad
\mathcal{D}_2=\sqrt 2(-1+\bold i)\partial_z,
\label{D1,2}
\end{eqnarray}
therefore the map
\begin{eqnarray}
\rho_{\mathbb C} (z)=
{1\over\sqrt 2(-1+\bold i)}\mathcal D_2
\quad ,\quad
\rho_{\mathbb C} (\overline{z})= 
{-1\over \sqrt 2(1+\bold i)}\mathcal D_1
\label{rho}
\end{eqnarray}
is nothing but the natural complex unitary Heisenberg algebra
representation, described in (15). Thus, by 
$
\mathcal D_1^*=\mathcal D_2\, ,\,
\mathcal D_2^*=\mathcal D_1 
$,
we have 

\begin{eqnarray}
\mathcal{PD}_Z^2={1\over 2}(\mathcal D_1\mathcal 
D_2,\mathcal D_2\mathcal D_1)=
H_Z-\lambda\sigma_0.
\label{PD^2}
\end{eqnarray}

This computation shows that the appearance of 
$\lambda\sigma_0$ in the squared operator
is due to Heisenberg's commutation relations. We also have 
$\mathcal{PD}_{Zf}^2=-{1\over 2}\square_\gamma -\lambda\sigma_0$, 
where the last term can be explained by the same argument.
Thus the square of the PaDi-operator is exactly the
classical Pauli operator.

{\bf (I2) The relativistic property of Pauli- and PaDi-operators.} 
The latter operator is characterized as the
non-relativistic spin operator, however, next we show that in
the situation given in this paper both the Pauli- and
the PaDi-operator are relativistic. This fact is proved below
by pointing out the exact matching of the PaDi-operator with the
original Dirac operator acting on 4-spinors.  

The complete Dirac operator, which involves also the potential $V$
due to the nucleus, is defined on a coordinate system $(t,x_1,x_2,x_3)$
on the 4-space by
\begin{eqnarray}
-{\hbar\over \bold ic}{\partial\over\partial t}
-eV -\sum_{r=1}^3\alpha_r
\big({\hbar\over\bold i}{\partial\over\partial x_r}+{e\over c}
\bold a^r\big)-
\alpha_0mc=-{\hbar\over \bold ic}{\partial\over\partial t}-H_D,
\end{eqnarray}
where $H_D$ is the Dirac-Hamilton operator, furthermore,
\begin{eqnarray}
   \alpha_1=
     \begin{pmatrix}
     0 & 0 & 0 & 1\\
     0 & 0 & 1 & 0\\
     0 & 1 & 0 & 0\\
     1 & 0 & 0 & 0
     \end{pmatrix} 
\quad , \quad
   \alpha_2=
     \begin{pmatrix}
     0 & 0 & 0 & -\bold i\\
     0 & 0 & \bold i & 0\\
     0 & -\bold i & 0 & 0\\
     \bold i & 0 & 0 & 0
     \end{pmatrix}\, ,\label{alpha1-4}\\ 
   \alpha_3=
     \begin{pmatrix}
     0 & 0 & 1 & 0\\
     0 & 0 & 0 & -1\\
     1 & 0 & 0 & 0\\
     0 & -1 & 0 & 0
     \end{pmatrix}
\quad ,\quad 
   \alpha_0=
     \begin{pmatrix}
     1 & 0 & 0 & 0\\
     0 & 1 & 0 & 0\\
     0 & 0 & -1 & 0\\
     0 & 0 & 0 & -1
     \end{pmatrix}\, .
\nonumber
\end{eqnarray}

The charged particle is considered in that unique inertial system where
the constant electromagnetic field defining the orbiting
spin has vanishing electric field and constant magnetic field
determined by the vector potential $\bold a=\lambda (-x_2,x_1,0)$.
Thus the relativistic 4-potential is $(0, \bold a^1,\bold a^2,0)$.
In this case the $H_D$ can be restricted onto the $(x_1,x_2)$-plane, 
meaning that the system is completely described by such 4-spinors which
are defined on the plane and for which both
the first- (corresponding to $t$) and the fourth-component 
(corresponding to $x_3$) vanish. Thus they are, actually, 
2-spinors defined by 
the second and third components. From the $4\times 4$-spin matrices 
only the $2\times 2$-matrices in the middle should be retained, since 
the rest part define only trivial operations. Note that these middle  
matrices are exactly the $\sigma$-matrices 
defined with respect to the perpendicular unit complex numbers 
$u_1=1$ and $u_2=-\bold i$.
Therefore, this restricted $H_D$ is nothing but 
the Pauli-Dirac operator $\mathcal{PD}$. 

The above Dirac equation clearly suggest the form of 
the PaDi-operator for bounded particles as well for free particles 
($V=0$).
In case of $H_Z$ the rest mass $m$ is neglected, while
for $H_{Zf}$ the rest-mass is defined by $m=\lambda /c$. 
  
This identification of the PaDi-operator with the Dirac operator proves 
the relativistic property of the PaDi-operator immediately. A general
Pauli operator or PaDi-operator can not be derived from the 
Dirac operator
in this simple way and they are indeed non-relativistic operators
which can be derived from the Dirac 
operator just by non-relativistic limit. 

Besides of finding the relativistic
version of the Schr\"odinger equation, Dirac's main interest 
was to find a linear equation which derives positive probability
densities. His main concern 
was that the second order relativistic
Klein-Gordon equation defined negative probabilities. 
It is noteworthy that Pauli was very
critical of this probabilistic argument of Dirac. According to him,
concept such as `` the probability of a particle to be at $\bold x$
in space" is meaningless for relativistic particles, thus it is
meaningless to interpret the wave $\xi(\bold x)$ as probability
amplitude. He regarded the Dirac equation, as well as the
Klein-Gordon equation, as the field equation
rather then as the equation of probability amplitude as Dirac 
preferred (cf. the famous Pauli-Weisskopf article where the authors 
pushed for the resurrection of the Klein-Gordon field by quantizing the 
KG-equation as well as the Maxwell equations). 

Pauli's concern about the relativistic probability amplitude is
solved in this article by identifying
the Dirac operator with the PaDi-operator. This operator lives
on the space determined
by the unique inertial system where $\bold E=0$ and $\bold B\not =0$
hold. Thus also the probability amplitude
is defined rather on the space then on the space-time. 

Also the hole-theory is completely avoidable by assuming
the existence of such particles only which are observable in the
mathematical model depicted here (observation practically means
performing the Stern Gerlach experiment).
Particles of negative energy have never been observed by this measuring
so far. Even transition from positive to negative energy state is 
meaningless, since the new particle is not observable by the
same observation even in the case when the new particle exist. 
For instance, to observe a negative charge
reduces the problem to choosing a complex structure $J$. Positive
charges are observable only on models defined by $-J$. Thus the result 
of a transition from a state to an other one of opposite charge is
not observable by the model defined by $J$. Yet, this argument does
not exclude the existence of positrons. On the contrary, these objects 
are most definitely predicted by this mathematical model.

{\bf (I3) Spectra of PaDi-operators; anomalous Zeeman zones.}
On $\mathcal S$, both the eigenfunctions and eigenvalues of 
$\mathcal{PD}$ can 
be explicitly determined. These computations should be carried out 
first for $\mathcal{PD}^2$. By representing
the $\bold C^2$-valued functions in the form  
$\phi =(\varphi_{1},\varphi_{2})$, the eigenfields of 
$-\square_\gamma -2\lambda \sigma_0$ 
with eigenvalue $\mu$ appear in the form 
$\phi_1=(\varphi ,0)$ or $\phi_2=(0,\varphi )$, where $\varphi$ is an 
eigenfunction of $-\square_\gamma$ with eigenvalue, say, $\nu$. Then 
the eigenvalue corresponding to $\phi_j$ is
$\mu_j =\nu +2(-1)^{j}\lambda$.  The explicit eigenvalues,
$-\lambda (4p+k+4\lambda k)$, 
of operator $\square_\gamma$ are known from (18).
Thus also the spectrum of $\mathcal{PD}^2$ is explicitly determined.

Then the eigenvalue problem regarding the
PaDi-operator $\mathcal{PD}_{Zf}$ can be easily completed. In fact,
all the eigenvalues $\mu_j$ above are strictly positive. Furthermore, 
for the fields
\begin{equation}\psi_{j+} =Q_j(\phi_j +
{1\over\sqrt\mu_j} \mathcal{PD}(\phi_j))\, ,\,
\psi_{j-} =Q_j(\phi_j -{1\over\sqrt\mu_j}\mathcal{PD}(\phi_j))
\end{equation}
we have 
\begin{equation}
\mathcal{PD}(\psi_{j+})=\sqrt\mu_j\psi_{j+}\, ,\, 
\mathcal{PD}(\psi_{j-})=-\sqrt\mu_j\psi_{j-}.
\end{equation}
(In these formulas the technical constant $Q_j$ is defined such that 
for a function satisfying
$||\varphi ||=1$ also $||\psi_{j\pm}||=1$ must be satisfied. 
Accordingly
$Q_j=1/((1-2(-1)^j|Z_\gamma |)^2+1)^{1\over 2}$.) Thus (60)
provides the eigenspinors of $\mathcal{PD}_{Zf}$ with the explicit
eigenvalues described in (61). 

The lowest eigenvalue, $2\lambda +4\lambda^2$, belongs to the
eigenfunctions $\varphi_{0}^{(a)}=\overline z^{(a)}$ on a zone 
$\mathcal H^{(a)}$. For $\mathcal{PD}^2_Z$ we have
$\mu_{10}=0$. This zero eigenvalue of the classical Pauli operator
is a current interest
in the literature. Note that, due to the additional constant 
$4\lambda^2$, this eigenvalue is non-zero for $\mathcal{PD}^2_{Zf}$. 
On microscopic level this additional term involves $\hbar^2$,
thus it is negligible for weak magnetic fields. It is not negligible,
however, for strong magnetic fields. 

Due to $\mu_{10}=0$, the eigenfunctions of $\mathcal{PD}_Z$
are well defined by (60) for all but this zero eigenvalue. To cover 
this missing case, the $\psi^{(a)}_{1+}$ is defined by 
$\phi_1$ and $\psi_{10-}^{(a)}=0$.
The constant $Q=1/\sqrt 2$ is independent of $j$.
Hereby, the eigenvalues and eigenfunctions both of 
$\mathcal{PD}_Z$ and $\mathcal{PD}_{Zf}$ 
are explicitly determined.

According to (60), there are two types of eigenspinors, $\psi_1$ and
$\psi_2$, depending on the position the generating function $\varphi$
is placed. The spinor spaces spanned by these eigenspinors are denoted
by $\mathcal S_1$ and $\mathcal S_2$ respectively. The fields in
these spaces represent the states of particles' position and
momentum respectively. The first question one should consider if
there is an overlapping between these two spaces? It turns out that 
$\mathcal S_2$ is completely contained in $\mathcal S_1$, furthermore,
$\mathcal S=\mathcal S_1$. Actually, these relations turn out to be true
in a more puzzling way: The eigenstates regarding the position and
momentum are described by the very same eigenspinors. 
One can use these relations to establish
the {\it uncertainty principle: The spinor fields can not be used at the
same time to describe both the position and momentum eigenstates. 
In other words, if one has a complete information about one type of
states, there is no information about the other type of states.} 

The anomalous zones $\mathcal S_\pm^{(a)}$, $\mathcal S^{(a)}_{1\pm}$, 
and $\mathcal S^{(a)}_{2\pm}$ are spanned by the appropriate 
eigenspinors derived from the zones $\mathcal H^{(a)}$ according to the
following formulas
\begin{eqnarray}
\mathcal S^{(a)}=
\mathcal H^{(a)}\times
\mathcal H^{(a)}=
\mathcal S^{(a)}_+ \oplus
\mathcal S^{(a)}_-=
\mathcal S^{(a)}_1= 
\mathcal S^{(a)}_{1+}\oplus 
\mathcal S^{(a)}_{1-},\\
\mathcal S_2^{(a)}= 
\mathcal H_{\mu_1 >0}^{(a)}\times\mathcal H^{(a)}=
\mathcal S^{(a)}_{2+}\oplus 
\mathcal S^{(a)}_{2-}\subset\mathcal S_1^{(a)}.
\end{eqnarray}
The whole spinor space is the direct sum of the anomalous zones, i. e.,
\begin{equation}
\mathcal S=\mathcal S_1=
\oplus_{a=0}^\infty\mathcal S^{(a)}=
\oplus_{a=0}^\infty\mathcal S_1^{(a)},\,\mathcal S_2=
\oplus_{a=0}^\infty\mathcal S_2^{(a)}\subset S_1.
\end{equation}

{\bf (I4) Anomalous kernels.}
The anomalous zones can be similarly analyzed than the normal ones. 
This theory includes explicit establishing
of anomalous kernels regarding projections, heat flows and PaDi-flows.
In this review only the projections onto the anomalous 
zone $\mathcal S^{(a)}_j$ is described. 
This $\bold C^2\otimes\bold C^2$-valued kernel
$
\mathcal Q^{(a)}_{(j)}(X,Y)
$
has the following component functions
\begin{eqnarray}
\mathcal Q^{(a)}_{(1)11}(X,Y)=
{\lambda\over 2\pi}(
L_a(\lambda |X-Y|^2)
e^{\lambda X\overline Y}+
(\lambda\overline XY)^a)
e^{-{\lambda\over 2} (|X|^2+|Y|^2)},
\\
\mathcal Q^{(a)}_{(1)22}(X,Y)=
{\lambda\over 2\pi}e^{\lambda X\overline Y}
L_a(\lambda |X-Y|^2)
e^{-{\lambda\over 2} (|X|^2+|Y|^2)},
\\
\mathcal Q^{(a)}_{(1)12}(X,Y)=
\mathcal Q^{(a)}_{(1)21}(X,Y)=0,\\
\mathcal Q^{(a)}_{(2)11}(X,Y)=
{\lambda\over 2\pi}(
L_a(\lambda |X-Y|^2)
e^{\lambda X\overline Y}
-(\lambda\overline XY)^a)
e^{-{\lambda\over 2} (|X|^2+|Y|^2)},\\
\mathcal Q^{(a)}_{(2)22}(X,Y)=
{\lambda\over 2\pi}
L_a(\lambda |X-Y|^2)
e^{\lambda X\overline Y}
e^{-{\lambda\over 2} (|X|^2+|Y|^2)}
,\\
\mathcal Q^{(a)}_{(2)12}(X,Y)=
\mathcal Q^{(a)}_{(2)21}(X,Y)=0.
\end{eqnarray}
By this kernel the concept of spinning point spreads can be introduced. 
More complicated formulas describe the projections $Q^{(a)}_{j\pm}$
onto the subzones $\mathcal S_{j\pm}^{(a)}$. Further complications
arise when the anomalous heat- and PaDi-flows are described.
Yet these computations are manageable which can be used for
establishing a well defined anomalous zonal Feynman measure
along with explicit stopwatch spinors. 

{\bf (J) Zonal Coulomb law, Lamb shift.}

The main reason for 
the Zeeman zones are established only for free particles ($V=0$) is
that they are not invariant with respect to multiplications 
with non-holomorphic functions such as 
the radial Coulomb potential 
$V=Q/r$.  
In subsequent papers the theory of zones will be extended
to bounded particles in two different ways. In one of them 
the quantum Coulomb 
operator $\bold V$ (multiplication with $V$) is modified such 
that it leaves the 
zones invariant. A natural modified operator is defined
by projecting $\bold V(\mathcal H^{(a)})$ 
back to $\mathcal H^{(a)}$. Then this zonal Coulomb operator,
$\bold V^{(a)}$, is an
integral operator with a smooth kernel which can be explicitly computed.
Note that the zonal Coulomb forth acts locally, thus the 
non-relativistic nature of the Coulomb-forth action is terminated.
Other problems caused by the original Coulomb law are terminated too.
For instance, there are eigenvalues with multiplicities $2$ (doublets)
in the spectrum of 
the Zeeman-Coulomb operator $H_Z+\bold V$. The existence of these 
doublets are argued in
the Lamb-Retherford experiment (Lamb shift). 
It is fascinating to see that this
multiplicity drops down to one with respect to each eigenvalue in the
global spectrum of the hydrogen atom having the zonal Coulomb potentials
$\bold V^{(a)}$. 
In other words, the real cause of the doublets in the
spectrum of the hydrogen atom is the non-relativistic Coulomb law.

There is an other way to build in the Coulomb law into the zonal theory,
namely, by building in the potential $V$ into the metric of a curved
Riemannian manifold. This idea is borrowed from Einstein's general 
relativity, obviously.
There is a natural generalization of  
2-step nilpotent Lie groups, leading to the concept of 
2-step nilpotent-type Riemannian manifolds, which can 
be used to carry out this idea.

A 2-step nilpotent Lie algebra is 
defined by an endomorphism
space, $E$, consisting of skew endomorphisms 
acting on a Euclidean space $\bold v$. 
In the place of the Euclidean space $\bold v$, consider   
a Riemannian manifold $(M,g)$ 
endowed by a smooth field, $E(p)$, of skew endomorphism 
spaces. If one assumes that the $E(p)$ is spanned by auto-parallel
complex structures, $J_Z(X)$, where $Z\in \bold z$, the system
$\{M,g,E(p)=J_{\bold z}(p)\}$ is called K\"ahler complex.
Such systems were investigated
in the literature just for dimensionality's $dim(\bold z)=1,3$
which correspond to K\"ahler resp. hyperk\"ahler manifolds. 

Principal bundle $(M,M\times\bold z,\bold z)$ corresponds to 
$(\bold v,\bold v\times\bold z,\bold z)$ in case of
nilpotent groups. The most important objects on this bundle 
are the {\it gauge connections}
$\omega$ which, by definition, satisfy the structure equation 
$d\omega =\Omega =\Omega^\alpha e_\alpha$, where the curvature form
$\Omega$ is defined by means of the endomorphisms 
$J_{e_\alpha}$. The gauge metrics $g_\omega$ on the gauge-bundle is
defined by the horizontal 
lift of metric
$g$ on $M$ to the horizontal subspace.  

Up to gauge transforms, 
this metric is unique, 
defining a 2-step nilpotent type Riemannian
manifold. The generalized Zeeman manifolds 
are defined by factorizing $\bold z$. 
The Zeeman manifolds defined by 
K\"ahler complexes form an important subclass. Particularly
interesting are the Bergman-Zeeman manifolds 
constructed on the bounded domains 
of $\bold C^n=\bold R^{2n}$. 
Zonal analysis along with several spectral investigations
is the next step in developing this theory. The Laplacian
on these Riemannian manifolds appears as a generalized Zeeman operator
with a ``built in electric potential". Thus the generalized 
zonal theory takes the contributions attributed 
to electric potentials into account.


\begin{thebibliography}{GSWH}


\bibitem[AB]{ab}
Y. Aharanov and D. Bohm:
\newblock Significance of electromagnetic potentials in the 
quantum theory.
\newblock {\em Phys. Rev.}, 115:485--491, 1959.

\bibitem[AC]{ac}
Y. Aharanov and A. Casher:
\newblock Ground state of spin-$1/2$ charged particle in a 
two-dimensional magnetic field.
\newblock{\em Phys. Rev. A},
19:2461-2462, 1979.


\bibitem[A]{a}
O. Alvarez:
\newblock Theory of strings with boundaries: fluctuations, topology
and quantum geometry.
\newblock {\em Nucl. Phys.}, 216:125--184, 1983.

\bibitem[AHS]{ahs}
J. Avron, I. Herbst, B. Simon:
\newblock Schr\"odinger operators with magnetic fields, I.,
General interaction.
\newblock{\em Duke Math. J.},
45:847--883, 1978.


\bibitem[Bo]{bo}
D. Bohm:
\newblock Quantum Theory,
\newblock Dover, 1979.

\bibitem[BSB]{bsb}
S. Duplij, W. Siegel, J. Bagger (editors):
\newblock Concise Encyclopedia of Supersymmetry and Noncommutative
Structures in Mathematics and Physics (SUSY Encyclopedia).
\newblock Kluwer, 2003.

\bibitem[Ee]{ee}
J. Eells:
\newblock Random walk on the fundamental group.
\newblock{\em Proc. Symp. Pure Math.},  27:211--217, 1975.

\bibitem[F]{f}
R. P. Feynman:
\newblock QED.
\newblock Princeton Univ. Press, 1988.


\bibitem[FH]{fh}
R. P. Feynman, A. R. Hibbs:
\newblock Quantum Mechanics and Path Integrals.
\newblock McGraw-Hill, 1965.

\bibitem[GSW]{gsw}
M. B. Green, J. H. Schwarz, E. Witten:
\newblock Superstring Theory. 
\newblock Cambridge Univ. Press, 1999.

\bibitem[H]{h}
S. W. Hawking:
\newblock Zeta function regularization of path integrals in 
curved space time.
\newblock{\em Comm. Math. Phys.}, 55:133--148, 1977.

\bibitem[Ka]{ka}
A. Kaplan:
\newblock Riemannian nilmanifolds attached to Clifford modules.
\newblock{\em Geom. Dedicata}, 11:127--136, 1981.

\bibitem[LL]{ll}
L. D. Landau, E. M. Lifshitz:
\newblock Quantum Mechanics.
\newblock Pergamon Press LTD, 1958.

\bibitem[Me]{me}
A. Messiah:
\newblock Quantum Mechanics.
\newblock Dover Publ. INC, 1999.

\bibitem[M\"u]{mu}
W. M\"uller:
\newblock Relative zeta functions, relative determinants and
scattering theory.
\newblock{\em Comm. Math. Phys.}, 192:309--347, 1998.

\bibitem[N]{n}
E. Nelson:
\newblock Feynman integrals and the Schr\"odinger equation.
\newblock{\em J. Math. Phys.}, 5(3):332--345, 1963.

\bibitem[OPS]{ops}
B. Osgood, R. Philips, P. Sarnak:
\newblock Extremals of determinants of Laplacians.
\newblock{\em J. Func. Anal.}, 80:148--211, 1988.

\bibitem[Po]{po}
A. M. Polyakov:
\newblock Quantum geometry of bosonic and Fermionic strings.
\newblock{\em Phys. Lett. B}, 103:207--213, 1981.

\bibitem[RS]{rs}
M. Reed, B. Simon:
\newblock Methods of Modern Mathematical Physics, II.
\newblock Academic Press, London, 1978.

\bibitem[RSz]{rsz}
F. Riesz, B. Sz\"okefalvi-Nagy:
\newblock Functional Analysis.
\newblock Dover, 1990.

\bibitem[Sch]{sch}
S. S. Schweber:
\newblock QED and the Man Who Made It: Dyson, Feynman, Schwinger, and
Tomonaga.
\newblock Princeton Univ. Press, 1994.

\bibitem[Schw]{schw}
J. Schwinger (Editor):
\newblock Selected Papers on Quantum Electrodinamics.
\newblock Dover Publ. INC, 1958.

\bibitem[Si]{si}
B. Simon:
\newblock Functional Integration and Quantum Physics.
\newblock Academic Press, London, 1979.

\bibitem[Sz1]{sz1}
Z. I. Szab\'o:
\newblock Locally non-isometric yet super isospectral spaces.
\newblock{\em Geom. funct. anal. (GAFA)}, 9:185--214, 1999.

\bibitem[Sz2]{sz2}
Z. I. Szab\'o:
\newblock Isospectral pairs of metrics on 
balls, spheres, and other manifolds with different local geometries.
\newblock{\em Ann. of Math.}, 154:437--475, 2001.

\bibitem[Sz3]{sz3}
Z. I. Szab\'o:
\newblock A cornucopia of isospectral pairs of metrics on 
spheres with different local geometries.
\newblock{\em Ann. of Math.}, 161:343--395, 2005.

\bibitem[Sz4]{sz4}
Z. I. Szab\'o:
\newblock Correction and improvement added to ``Isospectral pairs 
of metrics on balls and 
spheres with different local geometries".
\newblock{\em DG/0510202 (submitted)}

\bibitem[Sz5]{sz5}
Z. I. Szab\'o:
\newblock Normal zones on Zeeman manifolds with trace class heat and
Feynman kernels and well defined zonal Feynman integrals.
\newblock submitted.

\bibitem[Sz6]{sz6}
Z. I. Szab\'o:
\newblock Pauli-Dirac operators and anomalous zones on Zeeman manifolds.
\newblock submitted.


\bibitem[Sze]{sze}
G. Szeg\"o:
\newblock Orthogonal Polynomials.
\newblock AMS Providence, Rhode Island, vol(23), 1939.

\bibitem[Ta]{ta}
M. E. Taylor:
\newblock Partial Differential Equations. 
\newblock Springer, 1996.

\bibitem[To]{to}
S. Tomonaga:
\newblock The Story of Spin. 
\newblock Univ. of Chicago Press, 1997.

\bibitem[Ton1]{ton1}
A. Tonomura, N. T. Matsuda, R. Suzuki, A. Fukuhara, N. Osakabe,
H. Umezaki, J. Endo, K. Shinagawa, Y. Sugita, and H. Fujiwara:
\newblock Observation of Aharanov-Bohm effect by electron holography. 
\newblock {\em Phys. Rev. Lett.}, 48:1443--1446, 1982.

\bibitem[Ton2]{ton2}
A. Tonomura, N. Osakabe, T. Matsuba, T. Kawasaki, J. Endo, S. Yano,
and H. Yamada:
\newblock Evidence for Aharanov-Bohm effect with magnetic field 
completely shielded from electron wave. 
\newblock {\em Phys. Rev. Lett.}, 56:792-795, 1986.

\bibitem[V]{v}
M. Veltman:
\newblock Facts and Mysteries in Elementary Particle Physics. 
\newblock World Scientific, 2003.

\bibitem[W]{w}
S. Weinberg:
\newblock The Quantum Theory of Fields I, II. 
\newblock Cambridge Univ. Press, 1995.

\bibitem[We]{we}
H. Weyl:
\newblock The Theory of Groups and Quantum Mechanics. 
\newblock Dover, 1950.




\end{thebibliography}
\end{document}